\documentclass[preprint,sort&compress,final,10pt]{elsarticle}

\usepackage{amsmath}
\usepackage{amssymb}
\usepackage{graphicx}
\usepackage{latexsym}
\usepackage{epstopdf}
\usepackage{natbib}
\usepackage{color}
\usepackage{hyperref}
\usepackage{amsmath,amscd}

\newtheorem{theorem}{Theorem}[section]

\newtheorem{definition}[theorem]{Definition}
\newtheorem{lemma}[theorem]{Lemma}

\newtheorem{remark}[theorem]{Remark}
\numberwithin{equation}{section}

\def\Proof{\noindent{\bf Proof.}~}
\def\qed{\hfill$\square$\smallskip}

\def\Im{\mathrm{Im}}

\journal{\empty}
\date{}

\begin{document}

\begin{frontmatter}

\title{Almost periodic solutions for an asymmetric oscillation }

\author[au1]{Peng Huang}

\address[au1]{School of Mathematics Sciences, Beijing Normal University, Beijing 100875, P.R. China.}

\ead[au1]{hp@mail.bnu.edu.cn}

\author[au1]{Xiong Li\footnote{Corresponding author. Partially supported by the NSFC (11571041) and the Fundamental Research Funds for the Central Universities.}}

\ead[au1]{xli@bnu.edu.cn}

\author[au2]{Bin Liu\footnote{Partially supported by the NSFC (11231001).}}

\address[au2]{School of Mathematical Sciences, Peking University, Beijing 100871, P.R. China.}

\ead[au3]{bliu@pku.edu.cn}

\begin{abstract}
In this paper we study the dynamical behaviour of the differential equation
\begin{equation*}
x''+ax^+ -bx^-=f(t),
\end{equation*}
where $x^+=\max\{x,0\}$,\ $x^-=\max\{-x,0\}$, $a$ and $b$ are two different positive constants, $f(t)$ is  a real analytic almost periodic function. For this purpose, firstly, we have to establish some variants of the invariant curve theorem of planar almost periodic mappings, which was proved recently by the authors (see \cite{Huang}).\ Then we will discuss the existence of  almost periodic solutions and the boundedness of all solutions for the above asymmetric oscillation.
\end{abstract}


\end{frontmatter}

\section{Introduction}
In this paper we continue the work initiated in \cite{Ortega96} trying to understand  the dynamic behavior
 of an  asymmetric oscillation
\begin{equation}\label{a1}
\begin{array}{ll}
x''$+$a$$x^+$$- b$$x^-$=$f(t),
\end{array}
\end{equation}
where $a$,$b$ are two different positive constants,\ $x^+$=max$\{x,0\}$,\ $x^-$=max$\{-x,0\}$,\ $f(t)$ is a real analytic almost periodic function  with frequencies $\omega=(\cdots,\omega_\lambda,\cdots)_{\lambda\in \mathbb{Z}}\\ \in \mathbb{R}^\mathbb{Z}$, and $\omega=(\cdots,\omega_\lambda,\cdots)_{\lambda\in \mathbb{Z}}$ is a bilateral infinite sequence of rationally independent frequency, that is to say, any finite segments of $\omega$  are rationally independent. The general question that we have in mind is :\ under what conditions on $a$,\ $b$, \ $\omega$ and $f(t)$,\ (\ref{a1}) has almost periodic solutions and all solutions are bounded  ?

Due to the relevance with applied mechanics, for example,\ modeling some kind of suspension bridge\ (see \cite{Lazer90}), the following semilinear Duffing's equation
\begin{equation}\label{a2}
\begin{array}{ll}
$$x''$+$a$$x^+$$- b$$x^-$=$f(x,t)$$
\end{array}
\end{equation}
was widely studied,\ where $f(x,t)$ is a smooth $2\pi$-periodic function in $t$.

If the function $f(x,t)$ in (\ref{a2})\ depends only on the time $t$,\  equation\ (\ref{a2})\ becomes (\ref{a1}),
which had been studied by Dancer \cite{Dancer87} and Fu\v{c}ik \cite{Fucik80} in their investigations of boundary value problems associated to equations with jumping nonlinearities.\ For recent development,\ we refer to\ \cite{{Gallouet82}, {Habets93}, {Lazer89}}\ and references therein.

Ortega \cite{Ortega96} investigated the Lagrangian stability for the equation
\begin{equation}\label{a4}
x''+ax^+- bx^-=1+\gamma \,p(t),
\end{equation}
where  $\gamma$ is a small parameter.\ He proved that if  $|\gamma|$ is sufficiently small and $p \in \mathcal{C}^{4}(\mathbb{S}^{1})$ with $\mathbb{S}^1=\mathbb{R}/2\pi\mathbb{Z}$,  then all solutions are bounded,\ that is, for every solution $x(t)$,  it is defined for all $t\in \mathbb{R}$ and
$$\sup \limits_{t\in \mathbb{R} }(|x(t)|+|x{'}(t)|)<+\infty.$$
On the other hand,\ when
\begin{equation}\label{a5}
{1\over\sqrt a}+{1\over\sqrt b} \in \mathbb{Q},
\end{equation}
Alonso and Ortega \cite{Alonso96} constructed a $2\pi$-periodic function $f(t)$ such that all  solutions of  equation (\ref{a1}) with large initial conditions are unbounded.\ Moreover,\ for such a function $f(t)$,\  equation (\ref{a1}) has periodic solutions. This means that the unbounded solutions and periodic solutions coexist.

Liu \cite{Liu99} removed the smallness assumption on $|\gamma|$ in equation (\ref{a4}), and proved the boundedness of all solutions of equation (\ref{a1})  under  the resonant condition (\ref{a5})  and some other reasonable assumptions.

Liu \cite{Liu05} dealt with the existence of quasi-periodic solutions and the boundedness of all solutions of equation (\ref{a1}) when $f(t)$ is a real analytic,\  even and  quasi-periodic function with frequencies $\omega=(\omega_1,\omega_2,\cdots,\omega_n)$.\ Firstly,\ he obtained some invariant curve theorems of planar reversible mappings which are quasi-periodic in  the spatial variable.  As an application,\ he used the invariant curve theorem  to investigate the existence of  quasi-periodic solutions and the boundedness of all solutions of (\ref{a1}).

Recently,  in \cite{Huang2} we also considered the existence of quasi-periodic solutions and the boundedness of all solutions of equation (\ref{a1}) when $f(t)$ is a smooth  quasi-periodic function with the frequency $\omega=(\omega_1,\omega_2,\cdots,\omega_n)$.

In this paper we focus on equation (\ref{a1}) under the case that $f(t)$ is a real analytic almost periodic function  with frequencies $\omega=(\cdots,\omega_\lambda,\cdots)_{\lambda\in \mathbb{Z}}$. In order to obtain almost periodic solutions and the boundedness of all solutions of equation (\ref{a1}),  we first have to establish the invariant curve theorem for  almost  mappings.  Fortunately, recently we \cite{Huang} obtained  the  existence of  invariant curves of the following planar almost  mapping
\begin{equation}\label{M}
\mathfrak{M}:\quad \begin{array}{ll}
\left\{\begin{array}{ll}
x_1=x+y+f(x,y),\\[0.2cm]
y_1=y+g(x,y),
 \end{array}\right.\  (x,y)\in \mathbb{R} \times [a,b],
\end{array}
\end{equation}
where the perturbations $f(x,y)$ and $g(x,y)$ are almost periodic in $x$ with the frequency $\omega=(\cdots,\omega_\lambda,\cdots)_{\lambda\in \mathbb{Z}}$ and admit a rapidly converging Fourier series expansion, respectively.

However, the poincar\'{e} mapping of equation (\ref{a1}) does not have the form as $\mathfrak{M}$, but has the following expression
\begin{equation}\label{f11}
\begin{array}{ll}
\mathcal{M}_{\delta}:\ \ \left\{\begin{array}{ll}
x_1=x+\beta+\delta l(x,y)+\delta f(x,y,\delta),\\[0.2cm]
y_1=y+\delta m(x,y)+\delta g(x,y,\delta),
 \end{array}\right.\ \ \  (x,y)\in \mathbb{R} \times [a,b],
\end{array}
\end{equation}
where the functions $l,m, f, g$ are real analytic and  almost periodic in $x$ with the frequency $\omega=(\cdots,\omega_\lambda,\cdots)_{\lambda\in \mathbb{Z}}$, $f(x,y,0 ) =g(x,y, 0 ) = 0, $\ $\beta$ is a constant, $0<\delta< 1$ is a small parameter.

About the mapping $\mathcal{M}_{\delta}$ there are some results. When the functions $l,m, f, g$ in (\ref{f11}) are periodic and the mapping $\mathcal{M}_{\delta}$ possesses the intersection property, Ortega \cite{Ortega99}  obtained a variant of the small twist theorem and   also studied the existence of invariant curves of mappings with average small twist in  \cite{Ortega01}.

When the functions $l,m, f, g$ in (\ref{f11}) are quasi-periodic in $x$ and real analytic in a complex neighbourhood, the mappings $\mathcal{M}_{\delta}$ is  reversible with respect to the involution $\mathcal{G}:(x,y)\mapsto (-x,y)$,\ that is,\ $\mathcal{G} \mathcal{M}_{\delta} \mathcal{G} = \mathcal{M}_{\delta}^{-1}$. Liu  \cite{Liu05} established some variants of the invariant curve theorem which are similar to ones in \cite{Ortega99} and \cite{Ortega01}.

When the functions $l,m, f, g$ in (\ref{f11}) are quasi-periodic in $x$ and smooth functions, the mapping $\mathcal{M}_{\delta}$ possesses the intersection property, we \cite{Huang2} established some variants of the invariant curve theorem which  are also similar to ones in \cite{Ortega99} and \cite{Ortega01} in the smooth case, other than analytic case  on the basis of the  invariant curve theorem obtained in \cite{Huang1}.

Nevertheless the above results can not solve our problem. We have to establish the  existence of invariant curves of the  planar almost  mapping $\mathcal{M}_{\delta}$ on the basis of the  invariant curve theorem obtained in \cite{Huang}.

After we get those invariant curves theorems, as an application, we shall also use them to study the existence of almost periodic solutions and the boundedness of all solutions for an asymmetric oscillation (\ref{a1}).

Before ending  the introduction,\ we briefly recall some variants of the invariant curve theorem. When the perturbations in the mapping $\mathfrak{M}$ in (\ref{M}) is periodic and possesses the intersection property.\ Moser \cite{Moser62} considered  the twist mapping
$$
\mathfrak{M}_{0}:\quad \begin{array}{ll}
\left\{\begin{array}{ll}
x_1=x+\alpha(y)+\varphi_{1}(x,y),\\[0.2cm]
y_1=y+\varphi_{2}(x,y),
 \end{array}\right.
\end{array}
$$
where the functions $\varphi_{1},\varphi_{2}$ are assumed to be small and of periodic $2\pi$ in $x.$\ He obtained the existence of invariant closed curves of  $\mathfrak{M}_{0}$ which is of class $\mathcal{C}^{333}$.\ About  $\mathfrak{M}_{0}$,\ an analytic version of the invariant curve theorem was presented in \cite{Siegel97},\ a version in class $\mathcal{C}^{5}$ in  R\"{u}ssmann \cite{Russmann70} and a optimal version in class $\mathcal{C}^{p}$ with $ p>3$  in Herman \cite{Herman83,Herman86}.

When the perturbations $f(x,y),g(x,y)$ in (\ref{M}) are quasi-periodic in $x$, there are some results about the existence of invariant curves of the following planar quasi-periodic  mapping
\begin{equation}\label{a7}
\mathfrak{M}_{1}:\quad \begin{array}{ll}
\left\{\begin{array}{ll}
x_1=x+\alpha+y +f(x,y),\\[0.2cm]
y_1=y+g(x,y),
 \end{array}\right.\ \ \ \   (x,y)\in \mathbb{R} \times [a,b],
\end{array}
\end{equation}
where the functions $f(x,y)$ and $g(x,y)$ are quasi-periodic in $x$ with the frequency $\omega=(\omega_1,\omega_2$,$\cdots,\omega_n)$, real analytic in  $x$ and $y$, and $\alpha$ is a constant.

When the mapping $\mathfrak{M}_{1}$ in (\ref{a7}) is an exact symplectic map,\ $\omega_1,\omega_2,\cdots,\omega_n$, $2\pi\alpha^{-1}$ are sufficiently incommensurable, Zharnitsky \cite{Zharnitsky00} proved the existence of invariant curves of the map $\mathfrak{M}_{1}$.

When the mapping $\mathfrak{M}_{1}$ in (\ref{a7}) is  reversible with respect to the involution $\mathcal{G}:(x,y)\mapsto (-x,y)$,\ that is,\ $\mathcal{G} \mathfrak{M}_1 \mathcal{G} = \mathfrak{M}_1^{-1}, \omega_1,\omega_2,\cdots,\omega_n,2\pi\alpha^{-1}$ satisfy the Diophantine condition
\begin{equation*}
\Big|{\langle k,\omega \rangle {\alpha \over {2\pi}}-j}\Big|\geq {\gamma \over {|k|^\tau}},\ \ \ \ \forall\ \  k \in \mathbb{Z}^n\backslash\{0\},\ \ \forall j \in \mathbb{Z},
\end{equation*}
Liu  \cite{Liu05} established the invariant curve theorem for quasi-periodic reversible mapping $\mathfrak{M}_{1}$.

The rest of the paper is organized as follows.\ In Section 2,\ we  give some preliminaries. The nonresonant small twist theorem (Theorem \ref{thm6.4}) and the resonant small twist theorem (Theorem \ref{thm6.6}) for the almost periodic mapping $\mathcal{M}_{\delta}$  are proved in Sections 3. In Section 4,\ we will prove the existence of almost periodic solutions and the boundedness of all solutions for an asymmetric oscillation (\ref{a1}) depending almost periodically on time.

\section{Some preliminaries}

\subsection{The space of  real analytic almost periodic functions}
Our aim is to find almost periodic solutions $x$ for an  asymmetric oscillation (\ref{a1}) which admit a rapidly converging Fourier series expansion, thus we have to  define  the space of a kind of real analytic almost periodic functions.

We first define  the space of real analytic quasi-periodic functions $Q(\omega)$ as in \cite[chapter 3]{Siegel97}, here the $n$-dimensional frequency vector $\omega=(\omega_1,\omega_2,\cdots,\omega_n)$ is rational independent, that is, for any $k=(k_1,k_2,\cdots,k_n) \neq 0$,\ $\langle k,\omega \rangle =\sum k_j \omega_j \neq 0$.

\begin{definition}\label{def2.1}
A function $f:\mathbb{R} \rightarrow \mathbb{R}$ is called real analytic quasi-periodic with the frequency $\omega=(\omega_1,\omega_2,\cdots,\omega_n)$,  if there exists a real analytic function
$$F: \theta=(\theta_1,\theta_2,\cdots,\theta_n) \in \mathbb{R}^n \rightarrow \mathbb{R}$$
such that $f(t)=F(\omega_1t, \omega_2t,\cdots, \omega_nt)\ \text{for all}\ t\in \mathbb{R}$, where $F$ is $2\pi$-periodic in each variable and bounded in a complex neighborhood $\Pi_{r}^{n}=\{(\theta_{1},\theta_{2},\cdots,\theta_{n})\in \mathbb{C}^n : |\Im\ \theta_{j}|\leq r, j=1,2,\cdots, n \}$  of\, $\mathbb{R}^n$ for some $r > 0$. Here we call $F(\theta)$ the shell function of $f(t)$.
\end{definition}

Denote by $Q(\omega)$  the set of real analytic quasi-periodic functions with the frequency $\omega=(\omega_1,\omega_2,\cdots,\omega_n)$.  Given $f(t)\in Q(\omega)$, the shell function $F(\theta)$ of  $f(t)$ admits a Fourier series expansion
$$F(\theta)=\sum \limits_{k\in \mathbb{Z}^n} f_{k}e^{i \langle k,\theta \rangle },$$
where $k=(k_1,k_2,\cdots,k_n)$,  $k_j$ range over all integers and the coefficients $f_{k}$ decay exponentially with $|k|=|k_1|+|k_2|+\cdots+|k_n|$, then $f(t)$ can be represented as a Fourier series of the type from the definition,
$$
\begin{array}{ll}
f(t)=\sum \limits_{k\in \mathbb{Z}^n} f_{k}e^{i \langle k,\omega \rangle t}.
\end{array}
$$

In the following we define the norm of the real analytic quasi-periodic function $f(t)$ through that of the corresponding shell functions $F$.
\begin{definition}
For $r>0$, let $Q_{r}(\omega)\subseteq Q(\omega)$ be the set of real analytic quasi-periodic functions $f$ such that the corresponding shell functions $F$ are bounded on the subset $\Pi_{r}^{n}$\ with the supremum norm
$$\big|F\big|_{r}=\sup \limits_{\theta\in \Pi_{r}^{n}}|F(\theta)|=\sup \limits_{\theta\in \Pi_{r}^{n}}\Big|\sum_{k}f_{k}e^{i\langle k,\theta\rangle}\Big|<+\infty.$$
Thus we define $\big|f\big|_{r}:=\big|F\big|_{r}.$
\end{definition}

Similarly, one can give the definition of real analytic almost periodic functions with the frequency $\omega=(\cdots,\omega_\lambda,\cdots)_{\lambda\in \mathbb{Z}}$, which is a bilateral infinite sequence of rationally independent frequency, that is to say, any finite segments of $\omega$  are rationally independent. For this purpose, we first define analytic functions on some infinite dimensional space (see \cite{Dineen99}).
\begin{definition}
Let $X$ be a complex Banach space. A function $f : U\subseteq X \rightarrow \mathbb{C}$, where $U$ is an open subset of $X$, is called analytic if $f$ is continuous on $U$, and $f|_{U\cap X_1}$ is analytic in the classical sense as a function of several complex variables for each finite dimensional subspace $X_1$ of $X$.
\end{definition}

Note that for the bilateral infinite sequence of rationally independent frequency $\omega=(\cdots,\omega_\lambda,\cdots)_{\lambda\in \mathbb{Z}}$,
$$\langle k,\omega \rangle=\sum \limits_{\lambda\in \mathbb{Z}} k_\lambda \omega_\lambda,$$
where due to the spatial structure of the perturbation $k$ runs over integer vectors whose support
$$\mbox{supp}\, k=\big\{\lambda\ : \ k_{\lambda}\neq 0\big\}$$
is a finite set of $\mathbb{Z}$.

\begin{definition}\label{def2.2}
A function $f:\mathbb{R} \rightarrow \mathbb{R}$ is called real analytic almost periodic with the frequency $\omega=(\cdots,\omega_\lambda,\cdots)\in \mathbb{R}^\mathbb{Z}$, if there exists a real analytic function
$$F: \theta=(\cdots,\theta_\lambda,\cdots) \in \mathbb{R}^{\mathbb{Z}} \rightarrow \mathbb{R},$$
which admit a rapidly converging Fourier series expansion
$$F(\theta)=\sum \limits _{A\in\mathcal{S}} F_{A}(\theta),$$
where
$$F_{A}(\theta)={\sum \limits _{\mbox{supp}\, k \subseteq A}} f_{k}\,e^{i \langle k,\theta \rangle },$$
and $\mathcal{S}$ is a family of  finite subsets $A$ of $\mathbb{Z}$ with $\mathbb{Z}\subseteq \bigcup\limits_{A\in\mathcal{S}} A, \langle k,\theta \rangle=\sum \limits_{\lambda\in \mathbb{Z}} k_\lambda \theta_\lambda,$
such that $f(t)=F(\omega t)\ \text{for all}\ t\in \mathbb{R}$, where $F$ is $2\pi$-periodic in each variable and bounded in a complex neighborhood  $\Pi_{r}=\Big\{\theta=(\cdots,\theta_\lambda,\cdots) \in \mathbb{C}^{\mathbb{Z}} :\ |\Im\,\theta|_{\infty}\leq r\Big\}$ for some $r > 0$, where $|\Im\,\theta|_{\infty}=\sup \limits_{\lambda\in \mathbb{Z}} |\Im\,\theta_{\lambda}|$. Here $F(\theta)$ is called the shell function of $f(t)$.
\end{definition}

Suppose that the function $f(t)$ has the Fourier exponents $\{\Lambda_\lambda: \lambda\in \mathbb{Z}\}$, and its basis is $\{\omega_\lambda: \lambda\in \mathbb{Z}\}$. Then for any $\lambda\in \mathbb{Z}$, $\Lambda_\lambda$ can be expressed into
$$
\Lambda_\lambda=r_{\lambda_1}\omega_{\lambda_1}+\cdots+r_{\lambda_{j(\lambda)}}\omega_{\lambda_{j(\lambda)}},
$$
where $r_{\lambda_1},\cdots, r_{\lambda_{j(\lambda)}}$ are rational numbers. Therefore,
$$
\mathcal{S}=\{(\lambda_1,\cdots,\lambda_{j(\lambda)}):\lambda\in\mathbb{Z}\}.
$$
Thus, this family  $\mathcal{S}$ is not totally arbitrary.\ Rather,\ $\mathcal{S}$ has to be a spatial structure on $\mathbb{Z}$ characterized by the property that the union of any two sets in $\mathcal{S}$ is again in $\mathcal{S}$,\ if they intersect :
$$A,\ B \in \mathcal{S},\ \ \ A\cap B \neq \varnothing\ \ \ \Rightarrow \ \ \ A \cup B \in \mathcal{S}.$$
Moreover, we define
\begin{equation}\label{a00000}
{\mathbb{Z}}_{0}^{\mathbb{Z}}:=\Big\{k=(\cdots,k_\lambda,\cdots)\in\mathbb{Z}^{\mathbb{Z}} :\ \mbox{supp}\,k\subseteq A,\ A\in\mathcal{S}\Big\}.
\end{equation}

Denote by $AP(\omega)$  the set of real analytic almost periodic functions with the frequency $\omega$ which is given by Definition \ref{def2.2}. From the definition, given $f(t)\in AP(\omega)$, the shell function $F(\theta)$ admits a rapidly converging Fourier series expansion
$$F(\theta)=\sum \limits _{A\in\mathcal{S}} F_{A}(\theta),$$
where
$$F_{A}(\theta)={\sum \limits _{\mbox{supp}\, k \subseteq A}} f_{k}\,e^{i \langle k,\theta \rangle }.$$
From the definitions of the support $\mbox{supp}\, k$ and $A$, we know that $F_{A}(\theta)$ is real analytic and $2\pi$-periodic in $\big\{\theta_\lambda\ :\ \lambda\in A\big\}$.

As a consequence of the definitions of $\mathcal{S}$ and the support $\mbox{supp}\, k$ of $k$ and $\mathbb{Z}_{0}^{\mathbb{Z}}$,\ the Fourier series of the shell function $F(\theta)$ has another form
$$F(\theta)=\sum \limits _{k \in {{\mathbb{Z}}_{0}^{\mathbb{Z}} }} f_{k}e^{i \langle k,\theta \rangle },$$
then $f(t)$ can be represented as a Fourier series of the type
\begin{equation*}
f(t)={\sum \limits_{A\in \mathcal{S}}}\ {\sum \limits _{\mbox{supp}\, k \subseteq A}} f_{k}e^{i \langle k,\omega \rangle t},\quad\quad \langle k,\omega \rangle=\sum \limits_{\lambda\in \mathbb{Z}} k_\lambda \omega_\lambda
\end{equation*}
or
\begin{equation*}
f(t)=\sum \limits _{k \in {{\mathbb{Z}}_{0}^{\mathbb{Z}} }}f_{k}e^{i \langle k,\omega \rangle t}.
\end{equation*}
If we define
$$
f_{A}(t)={\sum \limits _{\mbox{supp}\, k \subseteq A}} f_{k}\,e^{i \langle k,\omega \rangle t},
$$
then
$$f(t)=\sum \limits _{A\in\mathcal{S}} f_{A}(t).$$
From the definitions of the support $\mbox{supp}\, k$ and $A$, we know that $f_{A}(t)$ is a real analytic quasi-periodic function with the frequency $\omega_A=\big\{\omega_\lambda\ :\ \lambda\in A\big\}$.

\subsection{The norms of real analytic almost periodic functions}

In the following we will give a kind of the norm for real analytic almost periodic functions. Before we describe the norm,  some more definitions  and  notations  are useful.

The main ingredient of our perturbation theory is a nonnegative weight function
$$[\,\cdot\,]\ :\ \ \ A\ \mapsto\ [A]$$
defined on $\mathcal{S}.$ The weight of a subset may reflect its
size,\ its location or something else.\ This, however,\ is immaterial for the perturbation theory itself.\ Here only the properties of monotonicity and subadditivity are required :
$$A\subseteq B \ \ \ \ \Rightarrow\ \ \ \  [A]\leq [B],$$
$$A\cap B \neq \varnothing \ \ \ \ \Rightarrow \ \ \ \ [A\cup B]+[A\cap B]\leq [A]+[B]$$
for all $A,\ B$ in $\mathcal{S}$.\ Throughout this paper, we always use the following weight function
$$[A]=1+\sum\limits_{i\in A}log^{\varrho}(1+|i|),$$
where $\varrho>2$ is a constant.

Now we can define the norm of the real analytic almost periodic function $f(t)$ through that of the corresponding shell function $F$ just like in the quasi-periodic case.

\begin{definition}
Let $AP_{r}(\omega)\subseteq AP(\omega)$ be the set of real analytic almost periodic functions $f$ such that the corresponding shell functions $F$ are bounded on the subset $\Pi_{r}$\ with the norm
$$\|F\|_{m,r}=\sum \limits _{A\in\mathcal{S}} |F_{A}|_{r}\, e^{m[A]}=\sum \limits _{A\in\mathcal{S}} |f_{A}|_{r}\, e^{m[A]}<+\infty,$$
where $m>0$ is a constant and
$$|F_{A}|_{r}=\sup\limits_{\theta\in\Pi_{r}}|F_{A}(\theta)|=\sup\limits_{\theta\in\Pi_{r}} \Bigg|\sum \limits _{\mbox{supp}\, k \subseteq A}f_{k}\, e^{i\langle k,\theta\rangle}\Bigg|=|f_{A}|_{r}.$$
Hence we define
$$\|f\|_{m,r}:=\|F\|_{m,r}.$$

If $f(\cdot,y)\in AP_r(\omega)$, and the corresponding shell functions $F(\theta,y)$ are real
analytic in the domain $D(r,s)=\big\{(\theta,y)\in \mathbb{C}^\mathbb{Z} \times \mathbb{C}\ :\ |\Im\,\theta|_{\infty}<r, |y-\alpha|<s\big\}$ with some $\alpha\in \mathbb{R}$, we define
$$\|f\|_{m,r,s}=\sum \limits _{A\in\mathcal{S}} |F_{A}|_{r,s}\, e^{m[A]}=\sum \limits _{A\in\mathcal{S}} |f_{A}|_{r,s}\, e^{m[A]},$$
where
$$|F_{A}|_{r,s}=\sup\limits_{(\theta,y)\in D(r,s)}|F_{A}(\theta,y)|=\sup\limits_{(\theta,y)\in D(r,s)} \Bigg|\sum \limits _{\mbox{supp}\, k \subseteq A}f_{k}(y)\, e^{i\langle k,\theta\rangle}\Bigg|=|f_{A}|_{r,s}.$$

\end{definition}

\subsection{The properties of  real analytic almost periodic functions}

In the following some properties of real analytic almost periodic functions are given.

\begin{lemma}\label{lem2.4}
The following statements are true:\\
$(i)$ Let $f(t),g(t)\in AP(\omega)$,\ then $f(t)\pm g(t), g(t+f(t))\in AP(\omega);$\\[0.1cm]
$(ii)$ Let $f(t)\in AP(\omega)$ and $\tau=\beta t +f(t)\ (\beta+f'>0)$,\ then the inverse relation is given by $t={\beta^{-1}}\tau +g(\tau)$  and $g\in AP({\omega / \beta})$.\ In particular,\ if $\beta=1$,\ then $g\in AP(\omega).$
\end{lemma}
The detail proofs of Lemma \ref{lem2.4} can be seen in \cite{Huang}, we omit it here.

\subsection{The frequencies of the perturbations }
In this paper, the frequency $\omega=(\cdots,\omega_\lambda,\cdots)$ of the perturbations is not only rationally independent with $|\omega|_{\infty}=\sup\{|\omega_\lambda|:\lambda\in \mathbb{Z}\}<+\infty$,  but also it satisfies the strongly nonresonant  condition.

In a crucial fashion the weight function determines the nonresonance conditions for the small divisors arising in this theory. As we will do later on by way of an appropriate norm,\ it   suffices to estimate these small divisors from below not only in terms of the norm of $k$
$$|k|=\sum \limits_{\lambda \in \mathbb{Z}} |k_{\lambda}|,$$
but also in terms of the weight of its support
$$[[k]]=\min\limits_{\mbox{supp}\, k \subseteq A \in \mathcal{S}}[A].$$
Then the nonresonance conditions read
\begin{equation}\label{b3}
\begin{array}{ll}
{{| \langle k,\omega \rangle|} \geq {c \over {\Delta([[k]])\Delta(|k|)}}},\ \ \ \ 0\neq k \in \mathbb{Z}_{0}^{{\mathbb{Z}}},
\end{array}
\end{equation}
where, as usual, $c$ is a positive parameter and $\Delta$ some fixed approximation function as described in the following. One and the same approximation function is taken here in both places for simplicity, since the generalization is straightforward. A nondecreasing function $\Delta\ :\ [1,\infty) \rightarrow [1,\infty)$ is called an approximation function,\ if
\begin{equation}\label{b15}
 {{\log\Delta(t)}\over t} \searrow 0,\ \ \ \ 1\leq t \rightarrow\infty,
\end{equation}
and
$$\int_{1}^{\infty} {{\log\Delta(t)}\over {t^2}}\, dt <\infty.$$
In addition,\ the normalization $\Delta(1)=1$ is imposed for definiteness.

In the following we will give a criterion for the existence of strongly nonresonant
frequencies. It is based on growth conditions on the distribution function
$$N_n(t)=\mbox{card}\ \big\{A\in\mathcal{S}\ :\ \mbox{card}(A)=n,\ [A]\leq t\big\}$$
for $n\geq 1$ and $t\geq 0$.

\begin{lemma}\label{lem2.11}
There exist a constant $N_0$ and an approximation function $\Phi$ such that
$$
\begin{array}{ll}
N_n(t)\leq\left\{\begin{array}{ll}
0,\ \ \ \ \ \ \ \ \ \  t<t_n\\[0.4cm]
N_0 \Phi(t),\ \ \ t\geq t_n
 \end{array}\right.
\end{array}
$$
with a sequence of real numbers $t_n$ satisfying
$$n\log^{\varrho-1} n\leq t_n\sim n\, \log^{\varrho} n$$
for $n$ large with some exponent $\varrho-1>1$, we say $a_n \sim b_n$, if there are two constants $c, C$ such that $c a_n\le b_n\le C a_n$ and $c, C$ are independent of $n$.
\end{lemma}
The detail proofs of Lemma \ref{lem2.11} can be seen in \cite{Huang}, we omit it here.

According to Lemma \ref{lem2.11}, there exist an approximation function $\Delta_0$ and a probability measure $\mu$ on the parameter space $\mathbb{R}^{\mathbb{Z}}$ with support at any prescribed point such that the  measure of the set of $\omega$ satisfying the following inequalities
\begin{equation}\label{b6}
\begin{array}{ll}
{{| \langle k,\omega \rangle|} \geq {c \over {\Delta_0([[k]])\Delta_0(|k|)}}},\ \ \ \ c>0,\ \forall\ 0\neq k \in \mathbb{Z}_{0}^{{\mathbb{Z}}}.
\end{array}
\end{equation}
is positive for a suitably small $c$,\ the proofs can be found in \cite{Poschel90},\ we omit it here.

Throughout this paper,\ we always assume that the frequency  $\omega=(\cdots,\omega_\lambda,\cdots)$  satisfies the nonresonance condition (\ref{b6}).

\subsection{Invariant curve theorem for the almost periodic  mapping  $\mathfrak{M}$}

In this subsection we formulate an useful invariant curve theorem for the almost periodic  mapping  $\mathfrak{M}$ given by (\ref{M}).

\begin{definition}\label{def2.40}
Let $\mathfrak{M}$ be a mapping given by (\ref{M}). It is said that $\mathfrak{M}$ has the intersection property if
$$\mathfrak{M}(\mathbf{\Gamma}) \cap \mathbf{\Gamma} \neq \emptyset$$
for every  curve $\mathbf{\Gamma}:x=\xi+\varphi(\xi),\ y=\psi(\xi)$,\ where $\varphi$ and $\psi$ are real analytic and almost periodic in $\xi$ with the frequency  $\omega=(\cdots,\omega_\lambda,\cdots)$.
\end{definition}

We choose a rotation number $\alpha$ satisfying the inequalities
\begin{equation}\label{b7}
\begin{array}{ll}
\left\{\begin{array}{ll}
a+\gamma\leq \alpha\leq b-\gamma,\\[0.4cm]
\Big|{\langle k,\omega \rangle {\alpha \over {2\pi}}-j}\Big|\geq {\gamma \over {\Delta([[k]])\Delta(|k|)}},\ \ \ \  \mbox{for all}\ \  k \in \mathbb{Z}_{0}^{{\mathbb{Z}}}\backslash\{0\},\ \  j \in \mathbb{Z}
 \end{array}\right.
\end{array}
\end{equation}
with some positive constant $\gamma$, where $\Delta$ is some approximation functions (see Theorem \ref{thm2.10}).

\begin{theorem}\label{thm2.10}
There is an approximation function $\Delta$ such that for suitable $\gamma$,\ the set of $\alpha$ satisfying (\ref{b7}) has positive measure.
\end{theorem}

\begin{theorem}\label{thm2.11}
Suppose that the  almost periodic mapping $\mathfrak{M}$ given by (\ref{M}) has the intersection property, and for every $y$, $f(\cdot,y),g(\cdot,y)\in AP_r(\omega)$ with $\omega$ satisfying the nonresonance condition (\ref{b6}), and the corresponding shell functions $F(\theta,y),G(\theta,y)$ are real
analytic in the domain $D(r,s)=\big\{(\theta,y)\in \mathbb{C}^\mathbb{Z} \times \mathbb{C}\ :\ |\Im\,\theta|_{\infty}<r, |y-\alpha|<s\big\}$ with $\alpha$ satisfying (\ref{b7}). Then for each positive $\bar{\varepsilon}$, there is a positive $\varepsilon_0=\varepsilon_0(\bar{\varepsilon},r,s,m,\gamma,\Delta)$ such that if $f,g$ satisfy  the following smallness condition
$$\|f\|_{m,r,s}+\|g\|_{m,r,s}<{{\varepsilon}}_0,$$
then the almost periodic mapping $\mathfrak{M}$ has an invariant curve $\mathbf{\Gamma _{0}}$ with the form
$$
\begin{array}{ll}
\left\{\begin{array}{ll}
x=x'+\varphi(x'),\\[0.2cm]
y=\psi(x'),
 \end{array}\right.\ \ \
\end{array}
$$
where $\varphi, \psi$ are almost periodic  with the frequency $\omega=(\cdots,\omega_\lambda,\cdots)$, and the invariant curve $\mathbf{\Gamma _{0}}$ is of the form $y=\phi(x)$ with $\phi\in AP_{r'}(\omega)$ for some $r'<r$, and $\|\phi-\alpha\|_{m',r'}<\bar{\varepsilon},\ 0<m'<m$. Moreover,\ the  restriction of $\mathfrak{M}$ onto $\mathbf{\Gamma _{0}}$ is
 $$\mathfrak{M}|_{\mathbf{\Gamma _{0}}}: x_{1}^\prime=x^\prime +\alpha.$$
\end{theorem}

\begin{definition}\label{def2.400}
Let $\mathfrak{M}$ be a mapping given by (\ref{M}). We say that $\mathfrak{M}:\mathbb{R} \times [a,b]\to \mathbb{R}^2$  is an exact symplectic if
$\mathfrak{M}$ is  symplectic with respect to the usual symplectic structure $dr\wedge d\theta$ and for every  curve $\mathbf{\Gamma}:\theta=\xi+\varphi(\xi),\ r=\psi(\xi)$,\ where $\varphi$ and $\psi$ are real analytic and almost periodic in $\xi$ with the frequency  $\omega=(\cdots,\omega_\lambda,\cdots)$, we have
$$
\lim_{T\to+\infty}\frac{1}{2T}\int_{-T}^T\, rd\theta=\lim_{T\to+\infty}\frac{1}{2T}\int_{-T}^T\, r_1d\theta_1.
$$
\end{definition}

\begin{lemma}\label{lem2.12}
If the mapping $\mathfrak{M}$ given by (\ref{M}) is an exact symplectic map, then it has intersection property.
\end{lemma}

The detail proofs of the Theorem \ref{thm2.10}, \ref{thm2.11} and Lemma \ref{lem2.12}  can be seen in \cite{Huang}, we omit it here.

\subsection{The small twist theorem}

In this subsection we formulate an useful  small twist theorem which is a variant of the invariant curve theorem (Theorem \ref{thm2.11}) for the almost periodic  mapping  $\mathfrak{M}$.

In many applications,\ one may meet the so called small twist mappings
\begin{equation}\label{f1}
\begin{array}{ll}
\mathfrak{M}_{\delta}:\ \ \left\{\begin{array}{ll}
x_1=x+\beta+\delta y+f(x,y,\delta),\\[0.2cm]
y_1=y+g(x,y,\delta),\\[0.1cm]
 \end{array}\right.\  (x,y)\in \mathbb{R} \times [a,b],
\end{array}
\end{equation}
where the functions $f$ and $g$ are real analytic functions, and  almost periodic in $x$ with the frequency $\omega=(\cdots,\omega_\lambda,\cdots)$,  $\beta$ is a constant, $0<\delta< 1$ is a small parameter.

We choose the number $\alpha$ satisfying the inequalities
\begin{equation}\label{f2}
\begin{array}{ll}
\left\{\begin{array}{ll}
a+\gamma\leq \alpha\leq b-\gamma,\\[0.4cm]
\Big|{\langle k,\omega \rangle {{\beta+\delta\alpha} \over {2\pi}}-j}\Big|\geq {\gamma \over {\Delta([[k]])\Delta(|k|)}},\ \ \ \  \mbox{for all}\ \  k \in \mathbb{Z}_{0}^{{\mathbb{Z}}}\backslash\{0\},\ \  j \in \mathbb{Z}
 \end{array}\right.
\end{array}
\end{equation}
with some positive constant $\gamma$, where $\Delta$ is some approximation functions.

\begin{theorem}\label{thm6.1}
Suppose that the  almost periodic mapping $\mathfrak{M}_{\delta}$ given by (\ref{f1}) has the intersection property, and  for every $y,\delta$, $f(\cdot,y,\delta),g(\cdot,y,\delta)\in AP_r(\omega)$ with $\omega$ satisfying the nonresonance condition (\ref{b6}), and the corresponding shell functions $F(\theta,y,\delta),G(\theta,y,\delta)$ are real
analytic in the domain $D(r,s)=\big\{(\theta,y)\in \mathbb{C}^\mathbb{Z} \times \mathbb{C}\ :\ |\Im\,\theta|_{\infty}<r, |y-\alpha|<s\big\}$ with $\alpha$ satisfying (\ref{f2}). Then for each positive $\bar{\varepsilon}$ there is a positive $\varepsilon_0=\varepsilon_0(\bar{\varepsilon},r,s,m,\gamma,\Delta)$ such that if $f,g$ satisfy  the following smallness condition
\begin{equation}\label{f3}
\|f(\cdot,\cdot,\delta)\|_{m,r,s}+\|g(\cdot,\cdot,\delta\|_{m,r,s}<\delta{{\varepsilon}}_0,
\end{equation}
then the almost periodic mapping $\mathfrak{M}_{\delta}$ has an invariant curve $\mathbf{\Gamma _{0}}$ with the form
$$
\begin{array}{ll}
\left\{\begin{array}{ll}
x=x'+\varphi(x'),\\[0.2cm]
y=\psi(x'),
 \end{array}\right.\ \ \
\end{array}
$$
where $\varphi, \psi$ are almost periodic  with frequencies $\omega=(\cdots,\omega_\lambda,\cdots)$, and the invariant curve $\mathbf{\Gamma _{0}}$ is of the form $y=\phi(x)$ with $\phi\in AP_{r'}(\omega)$ for some $r'<r$, and $\|\phi-\alpha\|_{m',r'}<\bar{\varepsilon},\ 0<m'<m$. Moreover,\ the  restriction of $\mathfrak{M}_{\delta}$ onto $\mathbf{\Gamma _{0}}$ is of the form
$$\mathfrak{M}_{\delta}|_{\mathbf{\Gamma _{0}}}:\ \ \ \  x_{1}^\prime=x^\prime+\beta+\delta\alpha.$$
\end{theorem}

\begin{remark}
If all the conditions of Theorem \ref{thm6.1} hold,  given any $\alpha$ satisfying the inequalities (\ref{f2}), there exists an invariant curve $\mathbf{\Gamma _{0}}$ of $\mathfrak{M}_{\delta}$, which is almost periodic  with the frequency  $\omega=(\cdots,\omega_\lambda,\cdots)$,\ and the restriction of $\mathfrak{M}_{\delta}$ onto $\mathbf{\Gamma _{0}}$ has the form
$$\mathfrak{M}_{\delta}|_{\mathbf{\Gamma _{0}}}:\   x_{1}^\prime=x^\prime +\beta+\delta\alpha .$$
\end{remark}

This is the so called small twist theorem.\ It is not a direct consequence of Theorem \ref{thm2.11},\ but one can use the same procedure in the proof of Theorem \ref{thm2.11} to prove it.\ Because there is nothing new in the proof,\ we omit it here.
\begin{remark}
The above conclusion is also true for the following mapping
\begin{equation}\label{f4}
\begin{array}{ll}
\left\{\begin{array}{ll}
x_1=x+\beta+\delta h(y)+f(x,y,\delta),\\[0.2cm]
y_1=y+g(x,y,\delta)
 \end{array}\right.
\end{array}
\end{equation}
with $h{'}(y)\neq 0$, if we change the conditions (\ref{f3}) into
$$
\begin{array}{ll}
M\Big(\|f(\cdot,\cdot,\delta)\|_{m,r,s}+\|g(\cdot,\cdot,\delta)\|_{m,r,s}\Big)<{\varepsilon}_0
\end{array}
$$
with $M=\max\Big\{|h{'}|_{s}\, , 1 \Big\}$.
\end{remark}

\section{Some variants of  the  small twist theorem}

After we get the small twist theorem, motivated by the results in \cite{{Ortega99}, {Ortega01}}, now we are going to consider the mapping
\begin{equation*}
\begin{array}{ll}
\mathcal{M}_{\delta}:\ \ \left\{\begin{array}{ll}
x_1=x+\beta+\delta l(x,y)+\delta f(x,y,\delta),\\[0.2cm]
y_1=y+\delta m(x,y)+\delta g(x,y,\delta),
 \end{array}\right.\ \ \  (x,y)\in \mathbb{R} \times [a,b],
\end{array}
\end{equation*}
where the functions $l,m, f, g$ are real analytic and  almost periodic in $x$ with the frequency $\omega=(\cdots,\omega_\lambda,\cdots)$, $f(x,y,0 ) =g(x,y, 0 ) = 0, $\ $\beta$ is a constant, $0<\delta< 1$ is a small parameter.

We shall be interested in the existence of invariant curves for the one-parameter family of almost periodic mappings  $\{ {\mathcal{M}_{\delta}} \}$.  Thus we assume that  $\mathcal{M}_{\delta}$ has the intersection property for each $\delta.$ Similar to the results in the above two references,
we can obtain the nonresoant and resonant small twist theorems.

\subsection{The nonresonant small twist theorem}

In this subsection we formulate the nonresonant small twist theorem (Theorem \ref{thm6.4}) for the almost periodic mapping $\mathcal{M}_{\delta}$. In other words, we will prove   the following theorem.

\begin{theorem}\label{thm6.4}
In the previous setting assume that for every  $A\in \mathcal{S},\omega_A=\big\{\omega_\lambda\ :\ \lambda\in A\big\}$ and $2\pi / \beta$ are rationally independent,\ and
$$\lim _{T \to \infty} {{1\over T} \int _0^T {\partial l\over \partial y}(x,y)dx} \neq 0.
$$
Then there exists $\delta_{0} > 0$ such that the mapping $\mathcal{M}_{\delta}$ has an invariant curve in the domain $\mathbb{R}\times[a,b]$ if
$0<\delta < \delta_{0}.$
The invariant curve is real analytic and almost periodic  with the frequency  $\omega=(\cdots,\omega_\lambda,\cdots)$.
\end{theorem}

\noindent\textbf{Proof}: The proof is similar to \cite[Theorem 1]{Ortega01}, we just give a sketch here.

Similar to \cite{Ortega01}, one can use the intersection property  to obtain
$$\lim \limits_{T\rightarrow\infty}{1 \over T}\int _{0}^{T}m(x,y)dx=0.$$
Expand $l$ and $m$ into Fourier series
$$l(x,y)={\sum \limits_{A\in \mathcal{S}}}\ {\sum \limits _{\mbox{supp}\, k \subseteq A}} l_{k}(y)e^{i \langle k,\omega \rangle x},\ \ \ \ m(x,y)={\sum \limits_{A\in \mathcal{S}}}\ {\sum \limits _{\mbox{supp}\, k \subseteq A}} m_{k}(y)e^{i \langle k,\omega \rangle x}.$$
From the assumption of Theorem \ref{thm6.4}, we have
$$l'_{0}(y)=\lim _{T \to \infty} {{1\over T} \int _0^T {\partial l\over \partial y}(x,y)dx} \neq 0.$$
Also the function $ \tilde{l}(x,y)=l(x,y)-l_{0 }(y) $ satisfies
$\lim \limits_{T\rightarrow\infty}{1 \over T}\int _{0}^{T} \tilde{l}(x,y)dx=0.$

Let
$$h_{1}(x,y)={\sum \limits_{A\in \mathcal{S}}\ {\sum \limits _{\substack{{\mbox{supp}\, k \subseteq A}\\ 0<\mu[[k]]+\rho|k|<N}}}} l_{k}(y)e^{i \langle k,\omega \rangle x},$$
$$h_{2}(x,y)={\sum \limits_{A\in \mathcal{S}}\ {\sum \limits _{\substack{{\mbox{supp}\, k \subseteq A}\\ 0<\mu[[k]]+\rho|k|<N}}}} m_{k}(y)e^{i \langle k,\omega \rangle x},$$
where $\mu,\rho$ are two small and $N$ is a large positive parameters.

It is well known that for any $\varepsilon{''}>0$,\ there exists a positive integer $N$ depending on $l,m,\mu,\rho$ and $\varepsilon{''}$ such that
\begin{align*}
&\|l(x,y)-h_{1}(x,y)-l_{0}(y)\|_{m-\mu,r-\rho,s}+\|m(x,y)-h_{2}(x,y)\|_{m-\mu,r-\rho,s}\nonumber \\[0.2cm]
&\leq e^{-N}(\|l\|_{m,r,s}+\|m\|_{m,r,s}\Big)<\varepsilon{''}
\end{align*}
for $0<\mu<m,0<\rho<r.$ Moreover, the  estimate
\begin{eqnarray}\label{f12}
&\|l(x,y)-h_{1}(x,y)-l_{0}(y)\|_{m,r,s}+\|m(x,y)-h_{2}(x,y)\|_{m,r,s}<\varepsilon{''}
\end{eqnarray}
obviously holds.

Consider the difference equations
\begin{equation}\label{f13}
\begin{aligned}
\Phi(x+\beta,y)-\Phi(x,y)+h_{1}(x,y)=0,\\[0.2cm]
\Psi(x+\beta,y)-\Psi(x,y)+h_{2}(x,y)=0
\end{aligned}
\end{equation}
for unknown functions $\Phi$ and $\Psi$.  Since   for every  $A\in \mathcal{S},\omega_A=\big\{\omega_\lambda\ :\ \lambda\in A\big\}$ and $2\pi / \beta$ are rationally independent,\ it is easy to verify that
$$
\begin{aligned}
\Phi(x,y)=-{\sum \limits_{A\in \mathcal{S}}\ {\sum \limits _{\substack{{\mbox{supp}\, k \subseteq A}\\ 0<\mu[[k]]+\rho|k|<N}}}} {l_{k}(y)  \over {{e^{{i\langle k ,\omega \rangle}\beta}-1}}}e^{i\langle k ,\omega \rangle x},\\
\Psi(x,y)=-{\sum \limits_{A\in \mathcal{S}}\ {\sum \limits _{\substack{{\mbox{supp}\, k \subseteq A}\\ 0<\mu[[k]]+\rho|k|<N}}}}{ m_{k}(y)  \over {{e^{{i\langle k ,\omega \rangle}\beta}-1}}}e^{i\langle k ,\omega \rangle x}
\end{aligned}
$$
is a solution of (\ref{f13}). Moreover there is a positive constant $\zeta_{0}(N)$ such that
$$\|\Phi\|_{m,r,s}+\|\Psi\|_{m,r,s}\leq \zeta_{0} \Big(\|l\|_{m,r,s}+\|m\|_{m,r,s}\Big).$$

Let
\begin{eqnarray*}
R_{1}(x,y)&=&\Phi(x+\beta,y)-\Phi(x,y)+\tilde{l}(x,y)\\[0.2cm]
&=&\Phi(x+\beta,y)-\Phi(x,y)+l(x,y)-l_{0}(y)
\end{eqnarray*}
and
$$R_{2}(x,y)=\Psi(x+\beta,y)-\Psi(x,y)+m(x,y).$$
Then by (\ref{f12}), we get
$$\|R_{1}\|_{m,r,s}+\|R_{2}\|_{m,r,s}\leq {\varepsilon {''} }.$$

Define the change of variables $\mathcal{U}$ by
$$
\begin{aligned}
\theta = x+\delta \Phi(x,y),\\
r = y+\delta \Psi(x,y).
\end{aligned}
$$
Then the transformed mapping\ \ $\mathcal{U} \circ \mathcal{M}_{\delta} \circ \mathcal{U}^{-1}$ is of the form
$$
\begin{array}{ll}
\left\{\begin{array}{ll}
\theta_1=\theta+\beta+\delta l_{0}(r)+\delta \phi_{1}\circ \mathcal{U}^{-1}(\theta,r,\delta),\\[0.2cm]
r_1=r+\delta\phi_{2}\circ \mathcal{U}^{-1}(\theta,r,\delta),
 \end{array}\right.
\end{array}
$$
where
\begin{align*}
\phi_{1}(x,y,\delta)&=f(x,y,\delta)+R_{1}(x,y)+D_{1}(x,y,\delta)+l_{0}(y)-l_{0}(y+\delta \Psi(x,y)),\\
\phi_{2}(x,y,\delta)&=g(x,y,\delta)+R_{2}(x,y)+D_{2}(x,y,\delta)
\end{align*}
and
$$D_{1}(x,y,\delta)=\Phi(x_{1},y_{1})-\Phi(x+\beta,y),$$
$$D_{2}(x,y,\delta)=\Psi(x_{1},y_{1})-\Psi(x+\beta,y).$$

Since\ $l,m, f, g$ are real analytic and  almost periodic  with the frequency $\omega=(\cdots,\omega_\lambda,\cdots)$, by the definitions of $\phi_{1},\phi_{2}$ and Lemma \ref{lem2.4},\ we know $\phi_{1},\phi_{2}$ are real analytic and almost periodic  with the frequency $\omega=(\cdots,\omega_\lambda,\cdots)$.\ Hence $\phi_{1}\circ \mathcal{U}^{-1}$ and $\phi_{2}\circ \mathcal{U}^{-1}$ are real analytic and  almost periodic in $x$ with the frequency $\omega=(\cdots,\omega_\lambda,\cdots)$,\ which is guaranteed by Lemma \ref{lem2.4} and the definition of  $\mathcal{U}$.\ Moreover,\ similar to \cite{Ortega01}, there exists a constant\ $M_{2}>1$ such that
\begin{align*}
&\|\phi_{1}\circ \mathcal{U}^{-1}(\cdot,\cdot,\delta)\|_{m,r,s}+\|\phi_{2}\circ \mathcal{U}^{-2}(\cdot,\cdot,\delta)\|_{m,r,s}\\[0.2cm]
\leq & M_{2}\Big\{\|\phi_{1}(\cdot,\cdot,\delta)\|_{m,r,s}+\|\phi_{2}(\cdot,\cdot,\delta)\|_{m,r,s}\Big\}.
\end{align*}

By the assumption of Theorem \ref{thm6.4},\ we know
$$l'_{0}(r)=\lim _{T \to \infty} {{1\over T} \int _0^T {\partial l\over \partial r}d\theta} \neq 0.$$
Choose
$$N_{1}=\max\Big\{|l'_0|_{s}\, ,\ 1 \Big\}.$$
Let
\begin{align*}
\tilde{\delta}&={{{\varepsilon}_0}\over {N_{1}M_{2}}}.
\end{align*}

Since $f(\cdot,\cdot,0)=g(\cdot,\cdot,0)=0$,  then
$$\lim_{\delta\rightarrow 0^{+}}\Big\{\|f(\cdot,\cdot,\delta)\|_{m,r,s}+\|g(\cdot,\cdot,\delta)\|_{m,r,s}\Big\}=0.$$
Thus there is a $\delta_{1}>0$\ such that for $\delta\in (0,\delta_{1})$,
$$\|f(\cdot,\cdot,\delta)\|_{m,r,s}+\|g(\cdot,\cdot,\delta)\|_{m,r,s}<\varepsilon',\ \ \ 0<\varepsilon'={\tilde{\delta}\over 4}<1.$$

Similar to \cite{Ortega01},\ we also can get an estimate of the form
$$\|l_{0}(y)-l_{0}(r)\|_{m,r,s}<\bar{\Omega}(\delta),\ \ \ \lim_{\delta\rightarrow 0^{+}}\bar{\Omega}(\delta)=0,$$
where $\bar{\Omega}$ is an appropriate modulus of continuity that depends on $l_{0}$ and $\Psi$. Because $\lim \limits_{\delta\rightarrow 0^{+}}\bar{\Omega}(\delta)=0$,\ there is a $\delta_{2}>0$\ such that
$$\bar{\Omega}(\delta)<{\tilde{\delta}\over 4},\ \ \  \forall\ \delta\in (0,\delta_{2}).$$

Let
$$N_{2}=\|l\|_{m,r,s}+\|m\|_{m,r,s}.$$
From the definitions of $D_{1}$\ and $D_{2}$, similar to \cite{Ortega01}, there exists $M_{1}>0$,\ for any $\delta\in (0,\delta_{1})$, we have
\begin{align*}
& \|D_{1}(\cdot,\cdot,\delta)\|_{m,r,s}+\|D_{2}(\cdot,\cdot,\delta)\|_{m,r,s}\\[0.2cm]
\leq& M_{1}\delta\Big\{\|l\|_{m,r,s}+\|m\|_{m,r,s}+\|f(\cdot,\cdot,\delta)\|_{m,r,s}+\|g(\cdot,\cdot,\delta)\|_{m,r,s}\Big\}\\[0.2cm]
\leq &M_{1}(N_{2}+1)\delta.
\end{align*}

Choose
\begin{align*}
\varepsilon'' &={\tilde{\delta}\over 4},\\
\delta_{0} &=\min\Big\{\delta_{1},\delta_{2},{\tilde{\delta}\over {4M_{1}(N_{2}+1)}}\Big\}.
\end{align*}
Then for any $\delta\in (0,\delta_{0})$,\ we get
\begin{eqnarray*}
\|\phi_{1}(\cdot,\cdot,\delta)\|_{m,r,s}+\|\phi_{2}(\cdot,\cdot,\delta)\|_{m,r,s}\leq (\varepsilon'+\varepsilon''+M_{1}(N_{2}+1)\delta+{\tilde{\delta}\over 4})< \tilde{\delta},
\end{eqnarray*}
which imply that
$$N_{1}\Big(\|\phi_{1}\circ \mathcal{U}^{-1}(\cdot,\cdot,\delta)\|_{m,r,s}+\|\phi_{2}\circ \mathcal{U}^{-2}(\cdot,\cdot,\delta)\|_{m,r,s}\Big)<N_{1}M_{2}\tilde{\delta}.$$
Then if $\delta\in (0,\delta_{0})$, we have
$$
\begin{array}{ll}
N_{1}\Big(\|\phi_{1}\circ \mathcal{U}^{-1}(\cdot,\cdot,\delta)\|_{m,r,s}+\|\phi_{2}\circ \mathcal{U}^{-2}(\cdot,\cdot,\delta)\|_{m,r,s}\Big)<  {\varepsilon}_0.
\end{array}
$$

Therefore, for $0<\delta<\delta_{0}$, this mapping meets all assumptions of Theorem \ref{thm6.1},  thus the transformed mapping  $\mathcal{U} \circ \mathcal{M}_{\delta} \circ \mathcal{U}^{-1}$ has invariant curves. Undoing the change of variables we obtain the existence of invariant curves of $\mathcal{M}_{\delta}$.\qed

\begin{remark}\label{rem6.5}
Actually it follows from the proof of Theorem \ref{thm6.4} that if all conditions of Theorem \ref{thm6.4} hold,\ then the mapping $\mathcal{M}_{\delta}$ has many invariant curves $\mathbf{\Gamma _{0}}$, which can be labeled by the form $\mathcal{M}_{\delta}|_{\mathbf{\Gamma _{0}}} : x_{1}^\prime=x^\prime+\beta+\delta \alpha,\ \alpha\in \Big[l_{0}(a)+\gamma,\ l_{0}(b)-\gamma\Big]$ of the restriction of $\mathcal{M}_{\delta}$ onto $\mathbf{\Gamma _{0}}.$\ In fact,\ given any $\alpha\in \Big[l_{0}(a)+\gamma,\ l_{0}(b)-\gamma\Big]$  satisfying the nonresonance condition
$$\Big|{\langle k,\omega \rangle {{\beta+\delta\alpha} \over {2\pi}}-j}\Big|\geq {\gamma \over {\Delta([[k]])\Delta(|k|)}},\ \ \ \  \mbox{for all}\ \  k \in \mathbb{Z}_{0}^{\mathbb{Z}}\backslash\{0\},\ \  j \in \mathbb{Z},$$
there exists an invariant curve $\mathbf{\Gamma _{0}}$ which is real analytic and almost periodic  with the frequency $\omega=(\cdots,\omega_\lambda,\cdots)$,\ and the restriction of $\mathcal{M}_{\delta}$ onto $\mathbf{\Gamma _{0}}$ has the form
$$\mathcal{M}_{\delta}|_{\mathbf{\Gamma _{0}}} : x_{1}^\prime=x^\prime+\beta+\delta \alpha.$$
\end{remark}

\subsection{The resonant small twist theorem}

Now we will discuss the resonant case, which means that there exists some set $A\in \mathcal{S}$ such that $\omega_A=\big\{\omega_\lambda\ :\ \lambda\in A\big\}$ and $2\pi / \beta$ are rationally dependent.
Denote by $\mathbb{{S}}$ the set of all $A\in \mathcal{S}$ such that there exists the integer vector $k\in\mathbb{Z}_{0}^{{\mathbb{Z}}}$ such that $\mbox{supp} k\subseteq A, \langle k ,\omega \rangle \beta \in 2\pi\mathbb{Z}$.

Now the functions $l$ and $m$ in (\ref{f11}) can be represented in the form
\begin{eqnarray*}
l(x,y):=\widetilde{l}(x,y)+\overline{l}(x,y)={\sum \limits_{A\in {\mathcal{S}\backslash\mathbb{{S}}}}}\ {\sum \limits _{\mbox{supp}\, k \subseteq A}} l_{k}(y)e^{i \langle k,\omega \rangle x}+{\sum \limits_{A\in \mathbb{S}}}\ {\sum \limits _{\mbox{supp}\, k \subseteq A}} l_{k}(y)e^{i \langle k,\omega \rangle x},
\end{eqnarray*}
\begin{eqnarray*}
m(x,y):=\widetilde{m}(x,y)+\overline{m}(x,y)={\sum \limits_{A\in {\mathcal{S}\backslash\mathbb{{S}}}}}\ {\sum \limits _{\mbox{supp}\, k \subseteq A}} m_{k}(y)e^{i \langle k,\omega \rangle x}+{\sum \limits_{A\in \mathbb{S}}}\ {\sum \limits _{\mbox{supp}\, k \subseteq A}} m_{k}(y)e^{i \langle k,\omega \rangle x}.
\end{eqnarray*}
Note that
$$e^{i{ \langle k,\omega \rangle} \beta}-1 \neq 0 ,\ \ \ \  {\mbox{supp}\, k \subseteq A},\ A\in{\mathcal{S}\backslash\mathbb{{S}}},$$
and
$$\overline{l}(x+\beta,y)\equiv\overline{l}(x,y),\ \ \ \   \overline{m}(x+\beta,y)\equiv\overline{m}(x,y).$$

\begin{theorem}\label{thm6.6}
In the previous setting  \mbox{supp}ose that the function $\bar{l}$  satisfies
$$\overline{l}(x,y) >0,\ \ \ \ \ {\partial{\overline{l}(x,y)}\over \partial y} > 0.$$
Moreover we assume that there is a real analytic function $I(x,y)\equiv I(x+\beta,y)$ satisfying
\begin{equation}\label{f17}
{\partial I(x,y)\over \partial y} > 0,
\end{equation}
\begin{equation}\label{f18}
\overline{l}(x,y){{\partial I}\over \partial x}(x,y)+\overline{m}(x,y){{\partial I}\over \partial y}(x,y)\equiv 0,
\end{equation}
and two numbers $\widetilde{a}$ and $\widetilde{b}$ such that
$$a< \widetilde{a} < \widetilde{b} < b$$
and
\begin{equation}\label{f19}
I_{\max}(a)<I_{\min}(\widetilde{a}) \leq I_{\max}(\widetilde{a})<I_{\min}(\widetilde{b}) \leq I_{\max}(\widetilde{b})<I_{min}(b),
\end{equation}
where
$$I_{\min}(y):=\min _{x \in \mathbb{R}} I(x,y),\ \ \ \ I_{\max}(y):=\max _{x \in \mathbb{R}} I(x,y).$$\\[0.0001cm]
Then there exists $\delta_0 > 0$ such that if $0<\delta < \delta_0$, the mapping $\mathcal{M}_\delta$ has an invariant curve which is real analytic  almost periodic  with the frequency $\omega=(\cdots,\omega_\lambda,\cdots)$. The constant  $\delta_0$ does depend only on $a,\ b,\ \widetilde{a}, \, \widetilde{b},\ l(x,y), m(x,y)$ and $I(x,y)$.
\end{theorem}
\noindent\textbf{Proof}: The proof is similar to \cite[theorem 3.1]{Ortega99}, we just give a sketch here.

Consider the following difference equations
\begin{equation}\label{f20}
\begin{aligned}
\tilde{\Phi}(x+\beta,y)-\tilde{\Phi}(x,y)+{\tilde{l}^{N}}(x,y)=0,\\[0.2cm]
\tilde{\Psi}(x+\beta,y)-\tilde{\Psi}(x,y)+{\tilde{m}^{N}}(x,y)=0,
\end{aligned}
\end{equation}
where
$${\tilde{l}^{N}}(x,y)={\sum \limits_{A\in{\mathcal{S}\backslash\mathbb{S}}}}\ {\sum \limits _{\substack{{\mbox{supp}\, k \subseteq A}\\ 0<\mu[[k]]+\rho|k|<N}}} l_{k}(y)e^{i \langle k,\omega \rangle x},$$
$${\tilde{m}^{N}}(x,y)={\sum \limits_{A\in{\mathcal{S}\backslash\mathbb{S}}}}\ {\sum \limits _{\substack{{\mbox{supp}\, k \subseteq A}\\ 0<\mu[[k]]+\rho|k|<N}}} m_{k}(y)e^{i \langle k,\omega \rangle x}.$$
Since $\langle k, \omega \rangle \beta  \notin 2\pi\mathbb{Z}$ for $ {\mbox{supp}\, k \subseteq A},\ A\in{\mathcal{S}\backslash\mathbb{{S}}}$, we know from the proof of Theorem \ref{thm6.4} that
$$
\begin{aligned}
\tilde{\Phi}(x,y)=-{\sum \limits_{A\in{\mathcal{S}\backslash\mathbb{S}}}}\ {\sum \limits _{\substack{{\mbox{supp}\, k \subseteq A}\\ 0<\mu[[k]]+\rho|k|<N}}} {l_{k}(y) \over {{e^{{i\langle k ,\omega \rangle}\beta}-1}}}e^{i\langle k ,\omega \rangle x},\\
\tilde{\Psi}(x,y)=-{\sum \limits_{A\in{\mathcal{S}\backslash\mathbb{S}}}}\ {\sum \limits _{\substack{{\mbox{supp}\, k \subseteq A}\\ 0<\mu[[k]]+\rho|k|<N}}}  {m_{k}(y) \over {{e^{{i\langle k ,\omega \rangle}\beta}-1}}}e^{i\langle k ,\omega \rangle x}
\end{aligned}
$$
is a solution of (\ref{f20}).

As we did in the proof of Theorem \ref{thm6.4}, under the transformation $\mathcal{U}_{1}$
$$
\begin{array}{ll}
\left\{\begin{array}{ll}
\theta = x+\delta \tilde{\Phi}(x,y),\\
r = y+\delta \tilde{\Psi}(x,y),
 \end{array}\right.
\end{array}
$$
the transformed mapping $\mathcal{U}_{1} \circ \mathcal{M}_{\delta} \circ \mathcal{U}_{1}^{-1}$ has the form
$$
\begin{array}{ll}
\left\{\begin{array}{ll}
\theta_1=\theta+\beta+\delta \bar{l}(\theta,r)+\delta \tilde{f}(\theta,r,\delta),\\[0.2cm]
r_1=r+\delta \bar{m}(\theta,r)+\delta \tilde{g}(\theta,r,\delta),
 \end{array}\right.
\end{array}
$$
where the functions $\tilde{f},\tilde{g}$ are very small if $\delta$ is sufficiently small and $N$ is very large.\ Moreover, $\tilde{f}$ and  $\tilde{g}$  are almost periodic in $x$ with the frequency  $\omega=(\cdots,\omega_\lambda,\cdots)$,\ which is guaranteed by Lemma \ref{lem2.4} and the definition of  $\mathcal{U}_{1}$.

In the following, we will construct another transformation $\mathcal{U}_{2}$ such that the transformed mapping $\mathcal{U}_{2} \circ\mathcal{U}_{1} \circ \mathcal{M}_{\delta} \circ \mathcal{U}_{1}^{-1}\circ\mathcal{U}_{2}^{-1}$ takes the form

$$
\begin{array}{ll}
\left\{\begin{array}{ll}
\tau_1=\tau+\beta+\delta \Omega(\rho)+\delta S_{1}(\tau,\rho,\delta),\\[0.3cm]
\rho_1=\rho+\delta S_{2}(\tau,\rho,\delta).
 \end{array}\right.
\end{array}
$$

For each $\theta\in \mathbb{R}$ and $h \in \mathbb{R}$ with $I(\theta,a)\leq h \leq I(\theta,b)$,  denote by $R=R(\theta,h)$ the unique solution of
$$I(\theta,R)=h.$$
The implicit function theorem and (\ref{f17}) imply that $R$ is well defined. Moreover,\ $R$ is $\beta$-periodic in $\theta$ and satisfies
$$
R(\theta,I(\theta,r))=r, \ \ \ \ \ \text {for all} \ \ \ (\theta,r)\in \mathbb{A}:=\mathbb{R}\times[a,b].
$$
The domain of $R$ obviously  contains the strip $\{(\theta,h):\bar{I}(a) \leq h \leq \underline{I}(b)\}.$

Next define the periodic function
$$T:[\bar{I}(a),\underline{I}(b)] \rightarrow \mathbb{R},\ \ \ \ T(h)=\int _{0}^{\beta} {d\theta \over {\bar{l}(\theta,R(\theta,h))}}.$$
This function is positive and has a negative derivative, namely,
$$
T'(h)=-\int _{0}^{\beta} {1 \over {\bar{l}(\theta,R)^{2}}}{\partial \bar{l} \over \partial r }(\theta,R){\partial R \over \partial h }d\theta<0.
$$
(Notice that ${\partial Y \over \partial h }=({\partial I \over \partial y})^{-1}> 0$).

The frequency function is defined by
\begin{equation}\label{f27}
\Omega(h)={\beta \over T(h)} \ \ \  \text {for} \ \ \ h \in [\bar{I}(a),\underline{I}(b)].
\end{equation}
This is a positive function with positive derivative.

Denote by $\tilde{\mathbb{A}}$ the following domain
$$\tilde{\mathbb{A}}=\{(\theta,r):\theta  \in  \mathbb{R},\tilde{a}\leq r\leq \tilde{b}\}$$
and define the function $K:\tilde{\mathbb{A}}\rightarrow \mathbb{R}$ by
$$K(\theta,r)=\int _{0}^{\theta} {ds \over {\bar{l}(s,R(s,I(\theta,r)))}},$$
which is well defined from (\ref{f19}).\ Moreover,\ it follows that
$$
K(\theta+\beta,r)=K(\theta,r)+T(I(\theta,r))\ \ \  \text {for all } \ \ \ (\theta,r)\in \tilde{\mathbb{A}}.
$$
The derivatives of $K$ are given by

$${\partial K \over \partial \theta}(\theta,r)={1 \over {\bar{l}(\theta,r)}}-{\partial I \over \partial \theta}(\theta,r)\int _{0}^{\theta} {1 \over {\bar{l}(s,R)^{2}}}{\partial \bar{l} \over \partial r }(s,R){\partial R \over \partial h }(s,I)ds,$$

$${\partial K \over \partial r}(\theta,r)=-{\partial I \over \partial r}(\theta,r)\int _{0}^{\theta} {1 \over {\bar{l}(s,R)^{2}}}{\partial \bar{l} \over \partial r }(s,R){\partial R \over \partial h }(s,I)ds.$$
From (\ref{f18}) we obtain
\begin{equation}\label{f29}
\bar{l}(\theta,r) {\partial K \over \partial \theta}+\bar{m}(\theta,r) {\partial K \over \partial r}=1.
\end{equation}

Introduce the new variables $(\rho,\tau)$ as follows
\begin{equation}\label{f30}
\rho=I(\theta,r),\ \ \ \ \tau=\Omega(I(\theta,r))K(\theta,r).
\end{equation}
By (\ref{f18}),\ (\ref{f29}),\ we have

\begin{equation}\label{f31}
\bar{l}(\theta,r) {\partial \tau \over \partial \theta}+\bar{m}(\theta,r) {\partial\tau \over \partial r}= \Omega \circ  I,
\end{equation}
and
\begin{equation}\label{f32}
\tau(\theta+\beta,r)=\tau(\theta,r)+\beta,
\end{equation}
\begin{equation}\label{f33}
{\partial \tau \over \partial \theta}(\theta+\beta,r)={\partial \tau \over \partial \theta}(\theta,r),\ \
{\partial \tau \over \partial r}(\theta+\beta,r)={\partial \tau \over \partial r}(\theta,r).
\end{equation}

We can now define the mapping $\mathcal{U}_{2}$
\begin{eqnarray*}
\tilde{\mathbb{A}}\mapsto \mathbb{R}^2,\ \  (\theta, r) \mapsto (\tau(\theta,r),\rho(\theta,r)),
\end{eqnarray*}
where  $\tau,\rho$ are defined by (\ref{f30}). The periodicity of $I$  and (\ref{f32}) imply that $\mathcal{U}_{2}$ satisfies
\begin{equation}\label{f34}
\mathcal{U}_{2}(\theta+\beta,r)=\mathcal{U}_{2}(\theta,r)+(\beta,0)\ \ \ \text {for all } \ \ \ (\theta,r)\in \tilde{\mathbb{A}}.
\end{equation}

Moreover,\ similar to \cite{Ortega99}, $\mathcal{U}_{2}$ can be expressed in the form
\begin{eqnarray*}
\begin{array}{ll}
\mathcal{U}_{2} :\left\{\begin{array}{ll}
\tau = \theta+\mathcal{U}_{2}^{1}(\theta,r),\\[0.3cm]
\rho= r+\mathcal{U}_{2}^{2}(\theta,r).
 \end{array}\right.
\end{array}
\end{eqnarray*}
From (\ref{f34}) we deduce that $\mathcal{U}_{2}^{1},\mathcal{U}_{2}^{2}$ are $\beta$-periodic in $\theta$.

In the following, we will give an expression of  $\mathcal{U}_{2} \circ \mathcal{U}_{1} \circ \mathcal{M}_{\delta} \circ \mathcal{U}_{1}^{-1} \circ \mathcal{U}_{2}^{-1}$. From (\ref{f18}), (\ref{f31}) and (\ref{f32}) it follows that
\begin{eqnarray*}
\tau_{1}-\tau &=& \tau(\theta+\beta +\delta \bar{l}+\delta \tilde{f},r+\delta \bar{m}+\delta \tilde{g})-\tau(\theta,r)\\[0.0001cm]
&=& \beta +\tau(\theta+\delta \bar{l}+\delta \tilde{f},r+\delta \bar{m}+\delta \tilde{g})-\tau(\theta,r)\\[0.2cm]
&=& \beta +\delta\Big[ {\partial \tau \over \partial \theta} \bar{l}+ {\partial \tau \over \partial r}\bar{m}\Big]+\delta\Big[ {\partial \tau \over \partial \theta} \tilde{f}+ {\partial \tau \over \partial r}\tilde{g}\Big]+\delta O_{1}(\delta)\\[0.2cm]
&=& \beta +\delta \Omega(I(\theta,r))+\delta\Big[ {\partial \tau \over \partial \theta} \tilde{f}+ {\partial \tau \over \partial r}\tilde{g}\Big]+\delta O_1(\delta)
\end{eqnarray*}
and
\begin{eqnarray*}
\rho_{1}-\rho &=& I(\theta_{1},r_{1})-I(\theta,r)\\[0.3cm]
&=& I(\theta+\beta+\delta \bar{l}+\delta \tilde{f},r+\delta \bar{m}+\delta \tilde{g})-I(\theta,r)\\[0.3cm]
&=& \delta\Big[ {\partial I \over \partial \theta} \bar{l}+ {\partial I \over \partial r}\bar{m}\Big]+\delta\Big[ {\partial I \over \partial \theta} \tilde{f}+ {\partial I \over \partial r}\tilde{g}\Big]+\delta O_{2}(\delta)\\[0.3cm]
&=& \delta\Big[ {\partial I \over \partial \theta} \tilde{f}+ {\partial I \over \partial r}\tilde{g}\Big]+\delta O_{2}(\delta).
\end{eqnarray*}
Therefore the transformed mapping $\mathcal{U}_{2}\circ  \mathcal{U}_{1} \circ \mathcal{M}_{\delta} \circ \mathcal{U}_{1}^{-1} \circ  \mathcal{U}_{2}^{-1}$ of\, $\mathcal{U}_{1} \circ \mathcal{M}_{\delta} \circ \mathcal{U}_{1}^{-1}$ under the transformation $\mathcal{U}_{2}$ takes the form
$$
\begin{array}{ll}
\left\{\begin{array}{ll}
\tau_1=\tau+\beta+\delta \Omega(\rho)+\delta \psi_{1}\circ \mathcal{U}_{2}^{-1}(\tau,\rho,\delta),\\[0.2cm]
\rho_1=\rho+\delta \psi_{2}\circ\mathcal{U}_{2}^{-1}(\tau,\rho,\delta),
 \end{array}\right.
\end{array}
$$
where
$$ \psi_{1}= \Big[ {\partial \tau \over \partial \theta} \tilde{f}+ {\partial \tau \over \partial r}\tilde{g}\Big]+O_{1}(\delta),\ \ \psi_{2}=\Big[ {\partial I \over \partial \theta} \tilde{f}+ {\partial I \over \partial r}\tilde{g}\Big]+O_{2}(\delta),$$
and the remainder term $O_{1}(\delta)$ is composed by $\delta,\tilde{f},\tilde{g}$ and second-order derivatives of $\tau$,\ the remainder term $O_{2}(\delta)$ is composed by second-order derivatives of $I$  and $\delta, \tilde{f}, \tilde{g}$.\ Hence,\ the remainder terms $O_{1}(\delta), O_{2}(\delta)$ are real analytic and satisfy $O_{1}(\delta), O_{2}(\delta)\rightarrow 0$ as $\delta\rightarrow 0$.

Since\ $I$ is $\beta$-periodic in $\theta$,\  its second-order derivatives are also $\beta$-periodic in $\theta$.\ By (\ref{f33}),\ it follows that the second-order derivatives of $\tau$ are $\beta$-periodic in $\theta$. In the previous setting we know $\tilde{f}, \tilde{g}$ are real analytic almost periodic in $\theta$ with the frequency $\omega=(\cdots,\omega_\lambda,\cdots)$ and assume that there exists set $A\in \mathcal{S},\omega_A=\big\{\omega_\lambda\ :\ \lambda\in A\big\}$ and $2\pi / \beta$ are rationally dependent,\
by the definitions of $\psi_{1}, \psi_{2}$ and Lemma \ref{lem2.4}, $\psi_{1}, \psi_{2}$ are real analytic and almost periodic in $\theta$ with the frequency $\omega=(\cdots,\omega_\lambda,\cdots)$.\ Hence $\psi_{1}\circ \mathcal{U}^{-1}_{2}$ and $\psi_{2}\circ \mathcal{U}^{-1}_{2}$  are real analytic  almost periodic  in $\tau$ with the frequency  $\omega=(\cdots,\omega_\lambda,\cdots)$,\ which is guaranteed by Lemma \ref{lem2.4} and the definition of  $\mathcal{U}_{2}$.\ Moreover,\ similar to \cite{Ortega99},\ we have
$$\big\|\psi_{1}\circ \mathcal{U}_{2}^{-1}(\cdot,\cdot,\delta)\big\|+\big\|\psi_{2}\circ \mathcal{U}_{2}^{-1}(\cdot,\cdot,\delta)\big\|\leq k\Big(\big\|\psi_{1}(\cdot,\cdot,\delta)\big\|+\big\|\psi_{2}(\cdot,\cdot,\delta)\big\|\Big).$$

The definitions of $\Omega$ yields that
$$\Omega{'}(\rho) \neq 0.$$
Let
\begin{align*}
M_{3}=\max\Big\{|\Omega{'}(\rho)|_{s}\, , 1 \Big\},\quad \quad\tilde{\delta}={{{\varepsilon}_0}\over {k M_{3}}}.
\end{align*}
From the definitions of $\Omega,\ K ,\ I,\ \tilde{f} ,\ \tilde{g}$, then there exists a constant  $\Delta_{0}$ such that if $\delta\in (0,\delta_{0})$,\ we have
\begin{eqnarray*}
\|\psi_{1}(\cdot,\cdot,\delta)\|_{m,r,s}+\|\psi_{2}(\cdot,\cdot,\delta)\|_{m,r,s} \leq \tilde{\delta},
\end{eqnarray*}
which imply that
$$M_{3}\Big(\|\phi_{1}\circ \mathcal{U}^{-1}(\cdot,\cdot,\delta)\|_{m,r,s}+\|\phi_{2}\circ \mathcal{U}^{-2}(\cdot,\cdot,\delta)\|_{m,r,s}\Big)<kM_{3}\tilde{\delta}.$$
Then if $\delta\in (0,\delta_{0})$, we have
$$
\begin{array}{ll}
M_{3}\Big(\|\phi_{1}\circ \mathcal{U}^{-1}(\cdot,\cdot,\delta)\|_{m,r,s}+\|\phi_{2}\circ \mathcal{U}^{-2}(\cdot,\cdot,\delta)\|_{m,r,s}\Big)<{\varepsilon}_0.
\end{array}
$$

Hence,\ for $0<\delta<\delta_{0}$, this mapping meets all assumptions of Theorem \ref{thm6.1}, and then the transformed mapping  $\mathcal{U}_{2} \circ \mathcal{U}_{1} \circ \mathcal{M}_{\delta} \circ \mathcal{U}_{1}^{-1} \circ \mathcal{U}_{2}^{-1}$ has invariant curves. Undoing the change of variables we obtain the existence of invariant curves of $\mathcal{M}_{\delta}$.\qed

One can obtain the same conclusions as Remark \ref{rem6.5}, we omit them here. Also, if
$$\overline{l}(x,y) <0,\ \ \ \ \ {\partial{\overline{l}(x,y)}\over \partial y} < 0,$$
and other conditions remain, then the statement of Theorem\ \ref{thm6.6} is also true.

\section{Application}

In this section we will apply the above results to an asymmetric oscillation.
Consider the following equation
\begin{equation}\label{g1}
\begin{array}{ll}
x''$+$a$$x^+$$- b$$x^-$=$f(t),
\end{array}
\end{equation}
where $a$, $b$ are two positive constants $(a \neq b)$,\ $x^+=\max\{x,0\}$,\ $x^-=\max\{-x,0\}$,\ $f\in AP(\omega)$ is real analytic almost periodic with the frequency $\omega=(\cdots,\omega_\lambda,\cdots)$ and admit a rapidly converging Fourier series expansion.

\subsection{Action and angle variables}

Introduce a new variable $y=x'$,\ then (\ref{g1}) is equivalent to the following planar system
\begin{equation}\label{g2}
\left\{\begin{array}{ll}
x'=y, \\[0.1cm]
y'=-ax^++bx^-+f(t).
\end{array}\right.
\end{equation}

Let $C(t)$ be the solution of the initial value problem
\begin{eqnarray*}
\left\{\begin{array}{ll}
x''+ax^{+}-bx^{-}=0, \\[0.1cm]
x(0)=1, \, x'(0)=0.
\end{array}\right.
\end{eqnarray*}
Then it is well known that $C(t) \in \mathcal{C}^2(\mathbb{R})$,  which can be given by
$$
C(t)=\left\{\begin{array}{ll}
 \cos {\sqrt a}t, &|t| \in [0,{\pi \over {2\sqrt a}}],\\[0.3cm]
-\sqrt{a \over  b} \sin {\sqrt b}\big(|t|-{\pi \over {2\sqrt a}}\big), &|t| \in [{\pi \over {2\sqrt a}},{{\pi \over {2\sqrt a}}+{\pi \over {2\sqrt b}}}].
\end{array}\right.
$$
Define $S(t)$ be the derivative of $C(t)$,\ then $S(t) \in \mathcal{C}^1(\mathbb{R}) $ and
\\(i) $C(-t)=C(t),S(-t)=-S(t)$.\\
(ii) $C(t)$ and $S(t)$ are $2\tilde{\omega} \pi$-periodic functions, $\tilde{\omega}={1\over 2}( {1 \over \sqrt a}+ {1 \over \sqrt b})$.
\\ (iii) $S^2(t)+ a(C^+(t))^2  + b(C^-(t))^2 \equiv a$.

For $r>0$,\ $\theta$\ (mod $2\pi$),\ we define the following generalized polar coordinates $T:(r,\theta)\rightarrow (x,y)$ as
\begin{eqnarray*}
\left\{\begin{array}{ll}
x=\varrho{r^{1\over 2}}C(\tilde{\omega} \theta), \\[0.1cm]
y=\varrho{r^{1\over 2}}S(\tilde{\omega} \theta),
\end{array}\right.
\end{eqnarray*}
where $\varrho=\sqrt{2\over{a\tilde{\omega}}}$. It is easy to check that $T$ is a symplectic transformation.

The Hamiltonian associated to (\ref{g2}) is expressed in cartesian coordinates as
$$H(t,x,y)={1\over 2}y^2+{a\over 2}(x^{+})^2+{b\over 2}(x^{-})^2-f(t)x.$$
In the new coordinates $(\theta, r)$, it becomes
$$H(t,\theta,r)={\tilde{\omega}}^{-1}r-\varrho{r^{1\over 2}}C(\tilde{\omega} \theta)f(t).$$
Thus, system (\ref{g2}) is changed into the following generalized polar coordinate system
\begin{equation}\label{g4}
\left\{\begin{array}{ll}
\theta '={\tilde{\omega}}^{-1}-{1\over 2}\varrho C(\tilde{\omega} \theta)f(t){r^{-{1\over 2}}}, \\[0.1cm]
r'=\tilde{\omega} \varrho S(\tilde{\omega} \theta)f(t){r^{1\over 2}}.
\end{array}\right.
\end{equation}

Since system (\ref{g4}) is not period in $t$, but is $2\pi$ period in $\theta$.  This fact will motivate us to interchange the role of $\theta$ and $t$. Now we change the role of the variable $t$ and $\theta$, and yields that
\begin{equation}\label{g5}
\left\{\begin{array}{ll}
{dt\over d\theta} =\big[{\tilde{\omega}}^{-1}-{1\over 2}\varrho C(\tilde{\omega} \theta)f(t){r^{-{1\over 2}}}\big]^{-1}, \\[0.2cm]
{dr\over d\theta} =\tilde{\omega} \varrho S(\tilde{\omega} \theta)f(t){r^{1\over 2}}\big[{\tilde{\omega}}^{-1}-{1\over 2}\varrho C(\tilde{\omega} \theta)f(t){r^{-{1\over 2}}}\big]^{-1},
\end{array}\right.
\end{equation}
this system is $2\pi$-periodic in the new time variable $\theta$. Let $r_{*}$ be a positive number such that
$${\tilde{\omega}}^{-1}-{1\over 2}\varrho {r_{*}^{-{1\over 2}}}\big|C\big|\big|f\big|>0.$$
System (\ref{g5}) is well defined for $r\geq r_{*}$. Let $\big(t(\theta),r(\theta)\big)$ be a solution of (\ref{g5}) defined in a certain interval $I=[\theta_{0},\theta_{1}]$ and such that $r(\theta)> r_{*}$ for all $\theta$ in $I$. The derivative ${dt\over d\theta}$ is positive and the function $t$ is a diffeomorphism from $I$ onto $J=[t_{0},t_{1}]$, where $t(\theta_{0})=t_{0}$ and $t(\theta_{1})=t_{1}$. The inverse function will be denoted by $\theta=\theta(t)$. It maps $J$ onto $I$.

\subsection{The expression of the Poincar\'{e} map of (\ref{g5})}

Our next goal is to obtain asymptotic expansions for $\theta_{1}$ and $r_{1}$. For $r_0$ large enough, the second equation of (\ref{g5}) can be rewritten as
$${d\over d\theta}r^{1\over 2}={1\over 2}{\tilde{\omega}}^{2}\varrho S(\tilde{\omega} \theta)f(t)+O(r^{-{1\over 2}}).$$
An integration of this equation leads to
\begin{equation}\label{g6}
\begin{array}{ll}
r(\theta)^{1\over 2}=r_0^{1\over 2}+O(1),& \theta \in [0,2\pi].
\end{array}
\end{equation}
Therefore,
$$r(\theta)^{-{1\over 2}}={{r_0}^{-{1\over 2}}}(1+{r_0}^{-{1\over 2}}O(1))^{-1}.$$
Expanding $(1+{r_0}^{-{1\over 2}}O(1))^{-1}$ yields
\begin{equation}\label{g7}
\begin{array}{ll}
{r(\theta)}^{-{1\over 2}}={r_0}^{-{1\over 2}}+O({r_0}^{-1}),& \theta \in [0,2\pi]
\end{array}
\end{equation}
for $r_0$ large enough,\ by (\ref{g7}) and the first equality of (\ref{g5}),\ we get
$${dt\over d\theta}={\tilde{\omega}}+O({r_0}^{-{1\over 2}}), \theta \in [0,2\pi],$$
which implies that
\begin{equation}\label{g8}
\begin{array}{ll}
t(\theta)=t_0+\tilde{\omega}\theta+O({r_0}^{-{1\over 2}}),& \theta \in [0,2\pi]
\end{array}
\end{equation}
for $r_0$ large enough.\ Substituting (\ref{g6}), \ (\ref{g7}),\ (\ref{g8}) into the second equality of (\ref{g5}),\ we have,\ for $\theta \in [0,2\pi]$,
\begin{eqnarray}\label{g9}
{d\over d\theta}r^{1\over 2}={1\over 2}{\tilde{\omega}}^{2}\varrho S(\tilde{\omega}\theta)f(t_0+\tilde{\omega}\theta)+O({{r_0}^{-{1\over 2}}})
\end{eqnarray}
for $r_0$ large enough.\ An integration of (\ref{g9}) over $\theta \in [0,2\pi]$ yields
$$
{r_1}^{1\over 2} =  {r_0}^{1\over 2}+ {1\over 2}{\tilde{\omega}}^{2}\varrho{\int _0^{2\pi}} f(t_0+\tilde{\omega}\theta)S(\tilde{\omega}\theta)d\theta+O({r_0}^{-{1\over 2}})
$$
for $r_0$ large enough, where $r_1=r(2\pi)$.\ Substituting (\ref{g7}),\ (\ref{g8}) into the first equality of (\ref{g5}),\ we have,\ for $\theta \in [0,2\pi]$
\begin{eqnarray}\label{g11}
{dt\over d\theta}=\tilde{\omega}+{1\over 2}{\tilde{\omega}}^{2}\varrho C(\tilde{\omega}\theta)f(t_0+\tilde{\omega}\theta){r_0}^{-{1\over 2}}+O({r_0}^{-1})
\end{eqnarray}
for $r_0$ large enough.\ An integration (\ref{g11}) over $\theta \in [0,2\pi]$ yields
$$
t_1=t_0+2\tilde{\omega} \pi+{1\over 2}{\tilde{\omega}}^{2}\varrho\ r_0^{-{1\over 2}}{\int _0^{2\pi}}{C(\tilde{\omega}\theta)f(t_0+\tilde{\omega}\theta)}d\theta+O(r_0^{-1}),
$$
where $t_1=t(2\pi)$.

Then the Poincar\'{e} map $P$ of (\ref{g5}) has the expansion
\begin{equation*}
P:
\left\{\begin{array}{ll}
t_1=t_0+2\tilde{\omega} \pi+{1\over 2}{\tilde{\omega}}^{2}\varrho\ r_0^{-{1\over 2}}{\int _0^{2\pi}}{C(\tilde{\omega}\theta)f(t_0+\tilde{\omega}\theta)}d\theta+O(r_0^{-1}), \\[0.2cm]
{r_1}^{1\over 2} =  {r_0}^{1\over 2}+ {1\over 2}{\tilde{\omega}}^{2}\varrho{\int _0^{2\pi}} f(t_0+\tilde{\omega}\theta)S(\tilde{\omega}\theta)d\theta+O({r_0}^{-{1\over 2}}).
\end{array} \right.
\end{equation*}

\subsection{Intersection property of the Poincar\'{e} map $P$ of (\ref{g5})}

Since the Poincar\'{e} map $P$ of (\ref{g5}) ia an exact symplectic map, similar to the proof of Lemma \ref{lem2.12}, it is not difficult to show that if $P$ is also monotonic twist map, that is,
\begin{equation}\label{g26}
L(t_0)={\int _0^{2\pi}}{f(t_0+\tilde{\omega}\theta)C(\tilde{\omega}\theta)}d\theta\neq 0,\ \ \ \text{for all}\ \ \  t_{0} \in \mathbb{R},
\end{equation}
then $P$ possesses the intersection property.  Thereinafter we always assume that $f(t)$ satisfies (\ref{g26}). In order to  prove this result,  we first give an useful lemma, and its proof is simple.

\begin{lemma}\label{lem 7.1}
If $r_0 = r(t_0)$ is real analytic almost periodic   in $t_0$  and $F(r_0,t_0)$ is real analytic almost periodic  in $t_0$ with the same frequencies, then $F(r(t_0),t_0)$ is also real analytic almost periodic in $t_0$ with the same frequencies.
\end{lemma}
\Proof We claim that the shell functions of  $r(t_0),F(r_0,t_0)$ are $R(\theta),\widetilde{F}(r_0,\theta)$ respectively, then the shell function of $F(r(t_0),t_0)$ is $\widetilde{F}(R(\theta),\theta).$ Indeed, since $r(t_0),F(r_0,t_0)\in AP(\omega)$, then we have
$$r(t_0)=\sum \limits _{k \in {{\mathbb{Z}}_{0}^{\mathbb{Z}} }}r_{k}e^{i \langle k,\omega \rangle t_0},\  R(\theta)=\sum \limits _{k \in {{\mathbb{Z}}_{0}^{\mathbb{Z}} }} r_{k}e^{i \langle k,\theta \rangle }$$
and
$$F(r_0,t_0)=\sum \limits _{k \in {{\mathbb{Z}}_{0}^{\mathbb{Z}} }} F_{k}(r_0)e^{i \langle k,\omega \rangle t_0},\ \widetilde{F}(r_0,\theta)=\sum \limits _{k \in {{\mathbb{Z}}_{0}^{\mathbb{Z}} }} F_{k}(r_0)e^{i \langle k,\theta \rangle }.$$
Hence,
\begin{eqnarray*}
F(r(t_0),t_0)=\sum \limits _{k \in {{\mathbb{Z}}_{0}^{\mathbb{Z}} }} F_{k}\Bigg(\sum \limits _{k \in {{\mathbb{Z}}_{0}^{\mathbb{Z}} }}r_{k}e^{i \langle k,\omega \rangle t_0}\Bigg)e^{i \langle k,\omega \rangle t_0}.
\end{eqnarray*}
From the definition of the shell function, we know that
\begin{align*}
\sum \limits _{k \in {{\mathbb{Z}}_{0}^{\mathbb{Z}} }} F_{k}\Bigg(\sum \limits _{k \in {{\mathbb{Z}}_{0}^{\mathbb{Z}} }}r_{k}e^{i \langle k,\theta \rangle }\Bigg)e^{i \langle k,\theta \rangle }=\sum \limits _{k \in {{\mathbb{Z}}_{0}^{\mathbb{Z}} }} F_{k}(R(\theta))e^{i \langle k,\theta \rangle}=\widetilde{F}(R(\theta),\theta)
\end{align*}
is the shell function of $F(r(t_0),t_0)$.\qed

Now we are going to prove the following lemma.

\begin{lemma}
If $f(t)$ satisfies (\ref{g26}), then Poincar\'{e} map $P$ of (\ref{g5}) has intersection property.
\end{lemma}
\Proof Since the mapping $P$ is an exact symplectic almost periodic monotonic twist map, according to the paper by Huang, Li and Liu \cite{Huang}, there is a function $H$ such that the mapping $P$ can be written by
\begin{equation*}
r_0 = -\frac{\partial }{\partial t_0} H(t_1-t_0, t_0),\quad r_1 = \frac{\partial}{\partial t_1}H(t_1-t_0,t_0),
\end{equation*}
where $H$ is real analytic almost periodic in the second variable.

Now we prove the intersection property of the mapping $P$, that is, given any continuous almost periodic curve $\Gamma : r_0=r(t_0)$,  we will prove that $P(\Gamma)\cap\Gamma\ne\emptyset.$  Define  two sets $\mathbb{B}$ and ${\mathbb{B}}_1$ : the set $\mathbb{B}$ is bounded by four curves $\big\{(t_0,r_0) : t_0=\theta\big\}$, $\big\{(t_0,r_0) : t_0=\Theta\big\}$, $\big\{(t_0,r_0) : r_0=r_{*}\big\}$ and $\big\{(t_0,r_0) : r_0=r(t_0)\big\}$, the set ${\mathbb{B}}_1$ is bounded by four curves   $\big\{(t_0,r_0) : t_0=\theta\big\}$, $\big\{(t_0,r_0) : t_0=\Theta\big\}$, $\big\{(t_0,r_0) : r_0=r_{*}\big\}$ and the image of $\Gamma$ under $P$. Here we choose $r_{*}<\min r(t_0)$ and $r_{*}$ is smaller than the image of $\Gamma$ under $P$ as $\theta\leq t_0\leq \Theta$. It is easy to show that the difference of the areas of ${\mathbb{B}}_1$ and $\mathbb{B}$ is
$$\Delta(\theta,\Theta)=\int _{\theta}^{\Theta}r_{1}d t_1-\int _{\theta}^{\Theta}r_0 d t_0=H(t_{1}(\Theta)-\Theta,\Theta)-H(t_{1}(\theta)-\theta,\theta).$$
From the definition of $P$ and the Lemma \ref{lem 7.1}, we know that $t_{1}(\Theta)-\Theta=f_1(r_0(\Theta),\Theta)$ is almost periodic in $\Theta$ and $t_{1}(\theta)-\theta=f_1(r_0(\theta),\theta)$ is almost periodic in $\theta$. Hence using Lemma \ref{lem 7.1} again, it follows that  $\Delta$ is almost periodic in $\Theta$ and $\theta$.

Hence there are at least two pairs of $(\theta_1,\Theta_1)$ and $(\theta_2,\Theta_2)$ such that $\Delta(\theta_1,\Theta_1)<0,\Delta(\theta_2,\Theta_2)>0,$ The intersection property of  $P$ follows from this fact.\qed

To apply the almost periodic invariant curve theorems obtained in Section 3, we need a fixed annulus. For this reason, introduce a new variable $v$ and a small parameter $\delta >0$ by
\begin{equation}\label{g13}
r={\delta^{-2}}v^{-2},\ \ \ v\in[1,2].
\end{equation}
Obviously,\ $r \gg 1 \Leftrightarrow \delta \ll 1$. Under this change of variable, the Poincar\'{e} map $P$ of (\ref{g5})  is transformed into
$$
P_{\delta}:
\left\{\begin{array}{ll}
t_1=t_0+2\tilde{\omega} \pi+\delta \Phi(t_0)v_0+\delta\widetilde{f}(t_0,v_0,\delta),& \\[0.2cm]
v_1=v_0+\delta\Psi(t_0){v_0^{2}}+\delta\widetilde{g}(t_0,v_0,\delta),  &
\end{array} \right.
$$
where
$$\Phi(t_0)={{\tilde{\omega}}^{2}\varrho\over 2}{\int _0^{2 \pi}}{f(t_0+\tilde{\omega}\theta)C(\tilde{\omega}\theta)}d\theta,$$
$$\Psi(t_0)= -{{\tilde{\omega}}^{2}\varrho\over 2}{\int _0^{2 \pi}} f(t_0+\tilde{\omega}\theta)S(\tilde{\omega}\theta)d\theta$$
and $ \widetilde{f},\widetilde{g}\rightarrow 0$ as $\delta\rightarrow 0.$

Since (\ref{g13}) is a homeomorphism for $r$ large enough, then $P_{\delta}$ also has the  intersection property when $\delta$ is sufficiently small. Now we apply Theorems \ref{thm6.4} and \ref{thm6.6} to prove  the existence of almost periodic solutions and the  boundedness of all solutions  for (\ref{g2}).

\subsection{The main results}

Assume that $f$ has the following Fourier series
$$f(t)={\sum \limits_{A\in \mathcal{S}}}\ {\sum \limits _{\mbox{supp}\, k \subseteq A}} f_{k}e^{i \langle k,\omega \rangle t}.$$

\begin{theorem}\label{thm7.1}
If for every finite subset $A$ of $\mathbb{Z},\ \omega_A=\big\{\omega_\lambda\ :\ \lambda\in A\big\}, 1 / \tilde{\omega}$ are rational independent and  $f(t)$ satisfies (\ref{g26}). Then system (\ref{g2}) has infinitely many almost periodic solutions and all solution are bounded.
\end{theorem}
\noindent\textbf{Proof}: If for every $A\in \mathcal{S},\omega_A=\big\{\omega_\lambda\ :\ \lambda\in A\big\}, 1 / \tilde{\omega}$ are rational independent, in this case from Theorem \ref{thm6.4}, we know that $P_{\delta}$ has invariant curves in the domain $\mathbb{R}\times [1,2]$ if $\delta$ is sufficiently small and
$$\lim _{T \to \infty} {1\over T}\int _0^T \Phi(t_0) d{t_0} \neq 0,$$
$$\lim _{T \to \infty} {1\over T}\int _0^T \Psi(t_0) d{t_0} = 0.$$
By (\ref{g26}), we know
$$\lim _{T \to \infty} {1\over T}\int _0^T \Phi(t_0) d{t_0} ={{\omega_0}^{2}\varrho\over 2}\lim _{T \to \infty}{1\over T} \int _0^T{\int _0^{2 \pi}}{f(t_0+\omega_0\theta)C(\omega_0\theta)}d\theta\, d t_0\neq 0$$
$\text{for all}\ \ \  t_{0} \in \mathbb{R}.$
Expand $f$ into the following Fourier series
$$f(t)=f_0 + {\sum \limits_{A\in \mathcal{S}}}\ {\sum \limits_{\substack{k\neq 0\\ \mbox{supp}\, k \subseteq A}}}  f_{k}e^{i \langle k,\omega \rangle t}.$$
By Fubini's theorem,\ it follows that
$$\lim _{T \to \infty} \int _0^T \Phi(t_0) d{t_0} =- {{f_{0}{\tilde{\omega}}\varrho}\over 2}\int _{0}^{2\tilde{\omega} \pi}S(\theta)d\theta =0.$$
Hence if the conditions of Theorem \ref{thm7.1} hold,\ then the existence of the invariant curves of Poincar\'{e} map $P_\delta$ is guaranteed by Theorem \ref{thm6.4},\ the invariant curves are  real analytic almost periodic  with frequencies $\omega=(\cdots,\omega_\lambda,\cdots)$.\ Undoing the change of variables we obtain the invariant curves of Poincar\'{e} map $P$.\ Then system (\ref{g2}) has infinitely many almost periodic solutions  as well as the boundedness of solutions.\qed

\begin{remark}\label{rem7.2}
It follows  from the proof of Theorem \ref{thm7.1} that if the conditions of Theorem \ref{thm7.1} hold,\ then system (\ref{g2}) has infinitely many almost periodic solutions with frequencies $\Big\{\omega=(\cdots,\omega_\lambda,\cdots), {{2\pi\tilde{\omega}}\over {{2\pi}\tilde{\omega}+\delta\alpha}}\Big\}$ satisfying the following nonresonance condition
\begin{align*}
&\Big|{\langle k,\omega \rangle {{{2\pi}\tilde{\omega}+\delta \alpha} \over {2\pi}}-j}\Big|\geq {\gamma \over {\Delta([[k]])\Delta(|k|)}},\ \ \ \  \mbox{for all}\ \  k \in \mathbb{Z}_{0}^{{\mathbb{Z}}}\backslash\{0\},\ \  j \in \mathbb{Z},\\[0.2cm]
&\alpha\in \Big[{{f_{0}{\tilde{\omega}}\varrho} {\sqrt a}}\big({1\over  a}-{1 \over  b}\big)+\gamma, 2{{f_{0}{\tilde{\omega}}\varrho} {\sqrt a}}\big({1\over  a}-{1 \over  b}\big)-\gamma\Big],
\end{align*}
and $\gamma, \delta$ are sufficiently small.
\end{remark}

\begin{theorem}\label{thm7.3}
If there exists finite subset $A$ of $\mathbb{Z},\ \omega_A=\big\{\omega_\lambda\ :\ \lambda\in A\big\}, 1 / \tilde{\omega}$ being rational dependent and  $f(t)$ satisfies (\ref{g26}). Denote by $\mathbb{S}$ the set of finite subset $A$ of $\mathbb{Z}$ such that the integer vectors $k\in\mathbb{Z}_{0}^{{\mathbb{Z}}}, \mbox{supp}\, k\subseteq A$ satisfying $\tilde{\omega}\sum \limits_{i\in A}{k_i\omega_i} \in \mathbb{Z}$ and $\tilde{\omega}\sum \limits_{i\in A}{k_i\omega_i} \notin \mathbb{Z}$ for $ A\notin\mathbb{{S}}, \mbox{supp}\, k\subseteq A, k\in\mathbb{Z}_{0}^{\mathbb{Z}}$, by $f_{\mathbb{S}}(t)$ the function
$$f_{ \mathbb{S}}(t)={\sum \limits_{A\in \mathbb{S}}}\ {\sum \limits _{\mbox{supp}\, k \subseteq A}} m_{k}(y)e^{i \langle k,\omega \rangle t}.$$
Moreover, if
$$\Phi_{\mathbb{S}}(t_0)={{\tilde{\omega}}^{2}\varrho\over 2}\int _0^{2\pi} f_{\mathbb{K}}(t_0+\tilde{\omega} \theta)C(\tilde{\omega}\theta)d\theta \neq 0,\ \ \  \text{for all}\ \ \ t_{0}\in \mathbb{R}.$$
Then system (\ref{g2}) has many almost periodic solutions and all solutions are bounded.
\end{theorem}
\noindent\textbf{Proof}:  In this case from Theorem \ref{thm6.6},\ $l=\Phi(t_0)v_0,\ m=\Psi(t_0){{v_0}^2}$ and
$$\overline{l}=\Phi_{\mathbb{S}}(t_0)v_0,$$
$$\overline{m}=\Psi_{\mathbb{S}}(t_0){{v_0}^2},$$
where
$$\Phi_{\mathbb{S}}(t_0)={{\tilde{\omega}}^{2}\varrho\over 2} \int _0^{2 \pi} f_{\mathbb{S}}(t_0+\tilde{\omega}\theta)C(\tilde{\omega}\theta)d\theta, $$
$$\Psi_{\mathbb{S}}(t_0)= -{{\tilde{\omega}}^{2}\varrho\over 2} \int _0^{2\pi} f_{\mathbb{S}}(t_0+\tilde{\omega}\theta)S(\tilde{\omega}\theta)d\theta. $$
Since
$$\Phi_{\mathbb{S}}(t_0)={{\tilde{\omega}}^{2}\varrho\over 2}\int _0^{2\pi} f_{\mathbb{S}}(t_0+\tilde{\omega}\theta)C(\tilde{\omega}\theta)d\theta \neq 0,\ \ \ \forall\  t_{0} \in \mathbb{R},$$
without loss of generality, we assume that
$$\Phi_{\mathbb{S}}(t_0)={{\tilde{\omega}}^{2}\varrho\over 2} \int _0^{2\pi} f_{\mathbb{S}}(t_0+\tilde{\omega} \theta)C(\tilde{\omega} \theta)d\theta > 0,$$
then
$$\overline{l}=\Phi_{\mathbb{S}}(t_0)v_0 >0,\ \ \ \ {\partial \overline{l} \over \partial v_0}=\Phi_{\mathbb{S}}(t_0)>0.$$
We may choose the function $I$ as
$$ I(t_0,v_0)= \exp{ \Big( {\int _0^{t_0}{\Psi_{\mathbb{S}}(t) \over \Phi_{\mathbb{S}}(t)} }dt \Big)}v_0.$$
Note that
$${\partial I\over \partial {t_0}}(t_0,v_0)=\exp{ \Big( {\int _0^{t_0}{\Psi_{\mathbb{S}}(t) \over \Phi_{\mathbb{S}}(t)} }dt \Big)} {\Psi_{\mathbb{S}}(t_0) \over \Phi_{\mathbb{S}}(t_0)}v_0,$$\\[0.0001cm]
$${\partial I\over \partial {v_0}}(t_0,v_0)=\exp{ \Big( {\int _0^{t_0}{\Psi_{\mathbb{S}}(t) \over \Phi_{\mathbb{S}}(t)} }dt \Big)}. $$
Therefore
$$\overline{l}(t_0,v_0){{\partial I}\over \partial {t_0}}(t_0,v_0)+\overline{m}(t_0,v_0){{\partial I}\over \partial {v_0}}(t_0,v_0)=0.$$\\[0.0001cm]
It is  not difficult to verify that all assumptions in Theorem \ref{thm6.6} are satisfied.\ Hence,\ the existence of the invariant curves of Poincar\'{e} map $P_\delta$ is guaranteed by Theorem \ref{thm6.6},\ which are real analytic almost periodic with frequencies  $\omega=(\cdots,\omega_\lambda,\cdots)$.\ Transformed back to the origin system this mean that the existence of the invariant curves of Poincar\'{e} map $P$.\ Then there are infinitely many almost periodic solutions  as well as the boundedness of all solutions of system (\ref{g2}).\qed

\begin{remark}
Similar to Remark \ref{rem7.2}, it follows  from the proof of Theorem \ref{thm7.3} that if the conditions of Theorem \ref{thm7.3} hold,\ then system (\ref{g2}) has infinitely many almost periodic solutions with frequencies $\Big\{\omega=(\cdots,\omega_\lambda,\cdots), {{2\pi\tilde{\omega}}\over {{2\pi}\tilde{\omega}+\delta\alpha}}\Big\}$  satisfying the following nonresonance condition
$$\Big|{\langle k,\omega \rangle {{{2\pi}\tilde{\omega}+\delta \alpha} \over {2\pi}}-j}\Big|\geq {\gamma \over {\Delta([[k]])\Delta(|k|)}},\ \ \ \  \mbox{for all}\ \  k \in \mathbb{Z}_{0}^{{\mathbb{Z}}}\backslash\{0\},\ \  j \in \mathbb{Z},$$
$$\alpha\in \Big[\Omega(1)+\gamma,\Omega(2)-\gamma\Big],$$
and $\Omega$ is defined by (\ref{f27}), $\gamma, \delta$ are sufficiently small..
\end{remark}

\bigskip

\section*{References}
\bibliographystyle{elsarticle-num}

\end{document}